\documentclass[opre]{informs3}

\usepackage{float}
\floatstyle{ruled}
\newfloat{algorithm}{htbp}{loa}
\floatname{algorithm}{Algorithm}
\usepackage[margin=1in]{geometry}
\let\theoremstyle\relax

\usepackage{amsmath,amssymb,amsthm,mathtools}
\usepackage{booktabs,array,longtable}

\usepackage{enumitem}
\usepackage{xcolor}
\usepackage{microtype}
\usepackage{capt-of}
 \usepackage[authoryear,round]{natbib}
\usepackage[colorlinks=true,linkcolor=blue!55!black,
            citecolor=blue!55!black,urlcolor=blue!55!black]{hyperref}
\newcommand{\E}{\mathbb E}
\newcommand{\1}{\mathbf 1}
\newcommand{\cE}{\mathcal E}
\newcommand{\cEzero}{\mathcal E_0}
\newcommand{\cS}{\mathcal S}
\newcommand{\cSI}{\mathcal S^I}
\newcommand{\cSA}{\mathcal S^A}
\newcommand{\KEEP}{\mathrm K}
\newcommand{\DROP}{\mathrm D}
\newcommand{\Ev}{\mathsf E}
\newtheorem{theorem}{Theorem}[section]
\newtheorem{proposition}[theorem]{Proposition}
\newtheorem{lemma}[theorem]{Lemma}
\newtheorem{corollary}[theorem]{Corollary}
\theoremstyle{definition}
\newtheorem{assumption}[theorem]{Assumption}
\numberwithin{equation}{section}
\newcommand{\Pp}{\mathbb P}

\newcommand{\deff}{d_{\mathrm{eff}}}

\usepackage{etoolbox}
\patchcmd{\paragraph}{\it}{\upshape\bfseries}{}{}
\theoremstyle{definition}

\theoremstyle{remark}
\newtheorem{remark}[theorem]{Remark}

\numberwithin{equation}{section}

\providecommand{\ECHowTheorems}{\setcounter{theorem}{0}}
\renewcommand{\ECHowEquations}{\setcounter{equation}{0}}

\title{TBD}
\author{TBD}
\date{\today}

\begin{document}

\RUNAUTHOR{TBD}

\RUNTITLE{Multi-Turn LLM Inference Optimization}

\TITLE{Routing and Cache Retention for Multi-Turn LLM Inference: Theory and Implementation}

\ARTICLEAUTHORS{
\AUTHOR{Sentao Miao\footnotemark[1]}
\AFF{Leeds School of Business, University of Colorado Boulder, Boulder, CO 80309, USA, \EMAIL{sentao.miao@colorado.edu}}

\AUTHOR{Tianjiao Nie\footnotemark[1]}
\AFF{Qiuzhen College, Tsinghua University, Beijing 100084, China, \EMAIL{nietj21@mails.tsinghua.edu.cn}}

 \AUTHOR{Yuan Zhou\footnotemark[1]\footnotemark[2]}
 \AFF{Yau Mathematical Sciences Center \& Department of Mathematical  Sciences, Tsinghua University, Beijing 100084, China; Beijing Institute of Mathematical Sciences and Applications, Beijing 101408, China, \EMAIL{yuan-zhou@tsinghua.edu.cn}}

 \AUTHOR{Song Yuan}
 \AFF{StepFun Intelligence, Shanghai 200232, China, \EMAIL{oliveryuan@stepfun.com}}

 \AUTHOR{Yuanwei Lu}
 \AFF{StepFun Intelligence, Shanghai 200232, China, \EMAIL{luyuanwei@stepfun.com}}

  \AUTHOR{Yibo Zhu}
 \AFF{StepFun Intelligence, Shanghai 200232, China, \EMAIL{zhuyibo@stepfun.com}}

\renewcommand{\thefootnote}{\fnsymbol{footnote}}
\footnotetext[1]{The co-first authors contributed equally to this work and are listed in alphabetical order.}
\footnotetext[2]{Corresponding author.}

} 
\ABSTRACT{Multi-turn AI agents execute iterative loops of large language model (LLM) inference and external tool calls, transforming cluster serving into a dynamic resource-allocation problem where routing and key--value (KV) cache retention jointly govern prefix reuse across inter-turn delays. We formulate a fluid approximation of this system and show that optimal admission, routing, and eviction are governed by endpoint-specific shadow prices for computation and memory together with session continuation values. Building on this structure, we propose an asymptotically optimal buffer-reserve policy that absorbs stochastic fluctuations using a vanishing share of capacity. To handle asynchronous arrivals, delayed load reports, and engine-local memory management in production clusters, we translate these principles into a practical controller named \textsc{PARK}, which pairs online projected subgradient price updates and one-turn migration lookahead with a local KV retention index ranked by expected reuse value per byte-epoch. In real GPU cluster experiments serving a 197-billion-parameter mixture-of-experts model on multi-turn software-engineering agent benchmarks (SWE and Terminal-Bench), \textsc{PARK} outperforms all five load-balancing and prefix-aware baselines across contrasting traffic regimes. On 16~H200 GPUs, \textsc{PARK} raises output throughput by $24.2\%$ over the strongest baseline ($336.9\%$ over round robin) on turn-intensive SWE trajectories and by $40.0\%$ over the strongest baseline ($63.3\%$ over round robin) on generation-intensive Terminal-Bench trajectories. Controlled single-replica ablations show that our retention index alone cuts prior-prefix token loss from $13.26\%$ to $2.17\%$ and increases throughput by $36.7\%$ over least-recently-used (LRU) eviction, while 128-GPU scale-up experiments sustain $50{,}789$ to $67{,}363$~output tokens/s across 16~replicas with prefix cache hit rates of $91.3\%$ to $97.4\%$.


}
\maketitle

\section{Introduction}
\label{sec:introduction}

Large language models (LLMs) support a growing range of activities,
including customer service, writing, software development, and
tool-using AI agents. Delivering these benefits requires
\emph{inference}: the computation performed by a trained model to
generate a response to a user's request. Whereas training develops a
model's capabilities once, inference puts those capabilities to use
each time a request is served, consuming processing capacity, memory,
and energy on graphics processing units (GPUs). Inference is therefore a
recurring operating requirement with substantial costs: OpenAI
reportedly spent \$8.4 billion on inference in 2025
\citep{muppidi2026inference}, and \citet{mckinsey2026inference}
projects that inference will account for about 60\% of AI data-center
demand by 2030. At this scale, modest reductions in the resources
required per request can yield substantial operating savings or allow
existing infrastructure to serve more demand. Realizing these gains
poses a natural operations research problem: how can a provider
complete more useful work with limited resources and uncertain demand?

We study a provider operating multiple inference servers, each of
which may contain several GPUs. Requests arrive sequentially, and the
provider must immediately decide whether to admit each request and
which server should process it. Admitted requests occupy computation
and memory resources while their responses are generated, sharing those
resources with unfinished work from earlier admissions. Memory is
particularly important because generating a response requires storing
intermediate computations whose size grows with the amount of text
being processed \citep{kwon2023pagedattention}.

This allocation problem is further complicated because many requests
belong to \emph{multi-turn} conversations, where each turn consists of a
prompt and its response. A later prompt builds on the conversation's
earlier text, creating an opportunity to reuse previous computation.
The saved intermediate computations are called the \emph{key--value (KV)
cache}. Keeping this cache reduces the work required for the next
turn, but occupies GPU memory while the user reads a response or
executes external tools. Therefore, the provider must also decide which idle
conversations deserve continued storage and which caches to evict,
meaning remove from memory. Eviction frees capacity for other work, at
the cost of rebuilding the lost computations if the conversation
resumes.

Crucially, request routing and cache eviction are closely coupled in multi-turn conversations. Serving a turn creates an updated cache at the chosen server, while eviction determines whether that cache remains available to reduce the cost of subsequent requests. Consider a conversation whose cache resides at a busy server while another server has spare capacity. Routing the next request to the cache-holding server avoids recomputation but adds to its load; routing it to the other server uses spare capacity and establishes an updated cache there, but requires rebuilding the prefix. Whether routing elsewhere is preferable depends on how long the destination’s congestion advantage persists and how likely the conversation is to continue. The eviction decision, in turn, determines which idle caches to retain by weighing the memory freed by eviction against the potential savings from future reuse. A retained cache occupies memory that could support other requests and saves computation only if the next request arrives before eviction and is routed to the server holding it.

Existing inference systems demonstrate the practical value of preserving caches across conversation turns \citep{gao2024cachedattention} and coordinating distributed routing to balance cache reuse and server load \citep{srivatsa2025preble}. Related optimization studies develop bid-price admission and routing in a single-turn model under resource constraints \citep{chen2026routing} and cache eviction policies for stochastic multi-turn conversations \citep{zhang2025tail}. These developments motivate a unified optimization and analysis of routing and cache retention. As the preceding discussion illustrates, the central challenge is to coordinate these decisions while accounting for how each shapes the costs and benefits of the other across conversation turns. We therefore ask:

\medskip
\begin{center}
\fbox{\parbox{6in}{\linespread{1.3}\selectfont
\centering
\itshape
How should a scheduler jointly optimize admission, routing, and idle-cache retention to maximize the throughput under server-specific computation and memory constraints?
}}
\end{center}

\medskip
We formulate this problem as a stochastic control model and develop a tractable fluid relaxation that links active service, idle waiting, and cache ownership across conversation turns. Resource shadow prices and session continuation values provide a common basis for coordinating these decisions, yielding structural characterizations of optimal fluid routing and retention. We further develop a fluid-based policy with an asymptotic performance guarantee and translate the theoretical principles into PARK, a practical controller for multi-turn inference. Our main contributions are as follows.

\begin{itemize}
\item \textbf{A fluid framework for jointly managing service and future
cache reuse.}
We formulate a stochastic control model for admission, routing, and
cache retention in multi-turn conversations that tracks each
session's service status, idle age, and cache location. Because
continuous relaxations of routing and retention fractions yield
nonlinear capacity constraints, we develop a fluid approximation
expressed through state--action flow volumes. This formulation yields a
tractable linear program while preserving the link between current
service and future cache availability.

\item \textbf{Structural characterization of routing and cache
eviction.}
The fluid solution provides a common economic basis for evaluating
routing and retention decisions: dual resource prices measure the cost
of computation and memory, while continuation values capture the
consequences for later turns. Their comparison explains when routing
should follow an existing cache and when rebuilding it elsewhere is
worthwhile. For retention, we derive sufficient conditions connecting
cache length, storage cost, and the conditional time until the next
prompt. These conditions show why intermediate-length caches can offer
a better trade-off than short caches (which are cheap to rebuild) or
very long caches (which occupy more memory than their future reuse
justifies), and characterize how idle age should govern eviction through the session's return prospect and rebuild cost, 
yielding a mechanism distinct from that of the traditional LRU policy.

\item \textbf{A fluid-based policy with an asymptotic performance
guarantee.}
We turn the fluid solution into an online buffer-reserve policy that
makes immediate routing decisions and periodically updates its
admission and cache-retention quotas while reserving capacity buffers
against stochastic fluctuations. Letting arrival rates and server
capacities grow proportionally with a scale parameter $n$, we prove an
$O(\sqrt{n\log n})$ shortfall in expected completed turns relative to
the fluid benchmark. Thus, the shortfall per unit of scale vanishes,
providing a performance guarantee in large stochastic systems.

\item \textbf{A practical controller (\textsc{PARK}) informed by the
theory and evaluated at scale.}
To handle asynchronous arrivals, delayed server snapshots, and
engine-local memory reclamation in production clusters, we translate our
routing and retention principles into an event-driven controller,
\textsc{PARK} (Price-based Admission and Routing with KV-cache
Retention). We evaluate \textsc{PARK} on two multi-turn
software-engineering agent workloads, serving a 197-billion-parameter
mixture-of-experts model on up to 128~H200 GPUs with a concurrency cap
of 4,096~sessions. In three-hour evaluations at this scale,
\textsc{PARK} sustains approximately 51,000--67,000~output tokens/s
across the two workloads. \textsc{PARK} outperforms standard routing
and caching baselines, while controlled single-replica experiments
show that its retention index mechanism reduces prior-prefix token loss from
13.26\% to 2.17\%.
\end{itemize}

The remainder of the paper is organized as follows. Section~\ref{sec:related-literature}
reviews the related literature. Section~\ref{sec:model} introduces the stochastic
model. Section~\ref{sec:fluid} develops the fluid formulation and structural
results. Section~\ref{sec:br-policy} presents the policy, its performance guarantee,
and its implementation. Section~\ref{sec:experiments} reports the experiments, and
Section~\ref{sec:conclusion} concludes. Proofs appear in the electronic companion.

\section{Related Literature}
\label{sec:related-literature}

Our work relates to three streams of research: LLM inference systems
and cache management, optimization for LLM inference, and dynamic
resource allocation.

\paragraph{LLM inference systems and cache management.}
Modern serving engines improve GPU memory efficiency and reuse
earlier computations across calls. PagedAttention reduces KV-cache
fragmentation and enables memory sharing \citep{kwon2023pagedattention},
SGLang reuses shared prompt prefixes across calls
\citep{zheng2024sglang}, CachedAttention preserves multi-turn KV
caches across a storage hierarchy \citep{gao2024cachedattention}, and
Preble coordinates distributed routing to balance cache locality and
server load \citep{srivatsa2025preble}. While these systems establish
the practical importance of cache reuse, we provide a unified
theoretical framework that jointly prices current resource use and
future conversation turns, explaining when routing and eviction
should depart from myopic heuristics.

Classical caching models also study retention decisions:
\citet{dehghan2016utility} optimize content-specific cache utilities,
and \citet{ferragut2016ttl} link optimal time-to-live timers to
inter-request time distributions. In multi-turn LLM inference,
however, idle caches compete directly with active requests for GPU
memory, and reuse value depends on whether and where future turns are
admitted. Our retention comparisons incorporate these service
decisions to characterize how cache length, idle age, and return
timing govern eviction.

\paragraph{Optimization and stochastic control for LLM inference.}
Recent work builds mathematical foundations for inference scheduling.
\citet{mitzenmacher2025queueing} survey connections among queueing,
service-time predictions, and LLM serving. \citet{dai2026throughput}
establish throughput optimality in batched processing networks,
including dependent AI-agent workflows. \citet{ao2026fluid} develop
fluid-guided admission and scheduling policies accounting for KV-memory
growth during generation, while \citet{lin2026heterogeneous} derive
asymptotically optimal admission and routing policies for competing
prompt and generation stages. Related studies examine online batching
under KV-memory limits \citep{jaillet2026scheduling}, non-clairvoyant
scheduling without advance output lengths
\citep{feng2026nonclairvoyant}, and cache-aware preemption for tail
latency \citep{li2026tailaware}.

Two closely related studies address routing and eviction separately:
\citet{chen2026routing} develop online bid-price routing in a single-turn model under
concurrent-request and KV-memory limits, and \citet{zhang2025tail}
propose a length-aware variant of least-recently-used eviction that
is optimal for a tail-latency objective under a stochastic model of
multi-turn conversations. \citet{wu2026randomization} model
routing--eviction interactions, but their routing algorithm greedily
minimizes predicted latency for the current request, and their
competitive guarantees apply to eviction rather than the combined
policy. In contrast, our fluid framework jointly optimizes admission, routing,
and retention across conversation turns: routing determines where
future caches reside, while retention governs the cost of future
routing choices. By accounting for these intertemporal effects,
our framework yields structural characterizations of the joint
decisions and an asymptotic optimality guarantee for our
buffer-reserve policy.

\paragraph{Dynamic resource allocation.}
Methodologically, our approach builds on bid-price control
\citep{talluri1998bidprice} and fluid reoptimization
\citep{jasin2012resolving,bumpensanti2020resolving}. Reusable-resource
models track capacity occupied for stochastic service durations and
released upon completion \citep{goyal2025reusable}. Multi-turn
inference adds a post-service state effect: completing a turn can
leave behind a cache that continues to consume memory while reducing
future execution work, both of which our routing and retention rules
price.

Our state-based LP also connects to weakly coupled Markov decision
processes, where local session dynamics share global capacity
constraints \citep{adelman2008relaxations,bertsimas2016decomposable}.
In this literature, \citet{brown2022signals} incorporate shared
external signals, and \citet{brown2025fluid} and
\citet{gast2024reoptimization} establish finite-horizon asymptotic
optimality guarantees for fluid reoptimization with heterogeneous
components. Building on these ideas, we link session states across
turns to connect resource shadow prices with cache-ownership values,
accommodating immediate routing, random service requirements, and
post-service retention decisions to achieve an $O(\sqrt{n\log n})$
shortfall from the fluid benchmark.

\section{Model and Decision Problem}\label{sec:model}

\subsection{Sessions, Endpoints, and Cache States}

We consider a discrete-time LLM inference system operating over a finite horizon of $H$ epochs, indexed by $t=0,\ldots,H-1$, with terminal time $H$. The system serves conversation sessions across a set of endpoints $\cE=\{1,\ldots,E\}$, where each endpoint represents a server housing one or more GPUs (processors used for LLM computation) with per-epoch execution capacity $B_e$ and memory capacity $G_e$. Upon arrival, each prompt must be immediately assigned to an endpoint or rejected; we write $\cEzero=\{0\}\cup\cE$, where action $0$ denotes rejection. After a rejection, the session may depart or remain idle without a cache, in which case it may submit a new prompt or retry the rejected one in a later epoch.

We call a submitted prompt a request. Each session consists of a sequence of turns, where a turn comprises an admitted request and its generated response. Each completed request therefore corresponds to one completed turn. Serving a turn involves processing the prompt and generating the response, which may span multiple epochs. A session is \emph{active} while its admitted turn is unfinished and \emph{idle} while waiting to submit another prompt. After a turn completes, its \emph{key--value (KV) cache} stores intermediate computations that can be reused by a later turn. Retaining an idle cache consumes endpoint memory between turns but reduces computation when the session returns. Each idle session has at most one retained cache: routing its next prompt to the endpoint holding that cache is a \emph{warm route}, whereas routing elsewhere, or routing without a retained cache, is a \emph{cold route} that requires recomputation.

Each live session has a state $s=(k,a,d,o,x)$ in a finite set $\cS$ (departed sessions are excluded), where:
\begin{itemize}
  \item $k$ records the session type or conversation phase, such as its context length, number of completed turns, or idle cache size;
  \item $a$ is the idle age (number of epochs elapsed since the previous turn completed), with $a=0$ for active sessions;
  \item $d$ is the attained service (number of epochs spent serving the unfinished turn), with $d=0$ for idle sessions;
  \item $o\in\cEzero$ is the endpoint holding an idle session's retained cache, with $o=0$ indicating no retained cache;
  \item $x\in\cEzero$ is the endpoint serving an unfinished turn, with $x=0$ for idle sessions.
\end{itemize}
Together, these coordinates capture all workload information needed to determine service requirements and return behavior. We partition $\cS=\cSI\mathbin{\dot\cup}\cSA$ into idle states ($\cSI$, where $d=x=0$) and active states ($\cSA$, where $a=0$); the cache-location coordinate $o$ is relevant only in $\cSI$.

Finally, let $M_s(t)$ be the number of sessions in state $s$ at the beginning of epoch $t$, and write $M(t)=(M_s(t))_{s\in\cS}$. When a turn completes at endpoint $e$, a continuing session becomes idle with its updated cache stored at $e$. Thus, the current routing decision directly determines where the next turn can reuse earlier computations.

\subsection{Decisions and Session Dynamics}

In each epoch $t$, decisions and state transitions occur in four steps:
\begin{enumerate}[label=(\arabic*)]
  \item \textbf{Observe the population.}
  The controller observes $M(t)$, including the serving endpoints of active sessions and the cache locations of idle sessions.

  \item \textbf{Route each arriving prompt.}
  Each beginning-of-epoch idle session in state $s$ independently submits a prompt with probability $p_s(t)$, while new sessions may enter with an initial prompt and no retained cache. Active sessions submit no prompts, so every session submits at most one prompt per epoch. Arrival times and new-session counts follow specified distributions, with new-session arrivals bounded per epoch.

  Indexing arrivals by $j=1,\ldots,J_t$ with session states $s_{t,j}$, the controller immediately chooses an action $e_{t,j}\in\cEzero$ upon each arrival, using all information observed up to that moment (including earlier arrivals and routing decisions in epoch $t$). Assigning $e_{t,j}\in\cE$ admits the prompt to endpoint $e_{t,j}$, where it remains until completion; action $0$ rejects it.

  \item \textbf{Observe service outcomes at the boundary.}
  Previously active and newly admitted turns receive service, yielding pre-retention states $\widetilde s_i(t+1)$ at the epoch boundary according to three state--action transition kernels:
  \begin{itemize}
    \item $Q_t^A(s,s')$ for an initially active session;
    \item $Q_t^I(s,s')$ for an idle session that submits no prompt;
    \item $P_t^e(s,s')$ for an arriving session assigned action $e\in\cEzero$.
  \end{itemize}
  Here $s'$ denotes the pre-retention state, and any probability mass unassigned to live states corresponds to session departure. Specifically, an unfinished turn remains active with attained service $d'=d+1$; a completed turn either departs or becomes idle at its serving endpoint $e$ with $(a',d',x',o')=(0,0,0,e)$; an unprompted idle session remains idle with age $a'=a+1$; and a rejected prompt releases any existing cache without starting service. Conditional on the beginning-of-epoch state and arrival sequence, transition and resource outcomes are independent across sessions and depend only on each session's current state and action, though computation, memory use, and completion may be correlated within a session.

  \item \textbf{Retain caches and update the population.}
  Let $\widetilde{\mathcal I}_e(t)$ be the set of idle sessions with a cache available at endpoint $e$ at the boundary (including newly completed turns). After observing service outcomes and cache sizes, the controller selects a subset $\mathcal I_e(t)\subseteq\widetilde{\mathcal I}_e(t)$ to retain and evicts the rest. Eviction clears only the cache location via $\Ev(k,a,0,o,0)=(k,a,0,0,0)$ without altering user conversation behavior. Post-retention states and the next population are then
  \begin{align*}
    s_i(t+1)
    &=\begin{cases}
      \Ev\bigl(\widetilde s_i(t+1)\bigr),
        &i\in\displaystyle\bigcup_{e\in\cE}
          \bigl(\widetilde{\mathcal I}_e(t)\setminus\mathcal I_e(t)\bigr),\\
      \widetilde s_i(t+1),&\text{otherwise},
    \end{cases}\\
    M_s(t+1)&=\sum_i\mathbf 1\{s_i(t+1)=s\},\qquad s\in\cS,
  \end{align*}
  where the sum ranges over all surviving and newly entering sessions.
\end{enumerate}

\subsection{Capacity Constraints and Optimization}

Let $\mathcal A_t$ be the set of active sessions at the beginning of epoch $t$. During epoch $t$, an initially active session $i\in\mathcal A_t$ incurs execution requirement $L_{t,i,e}$ at endpoint $e$ (with $L_{t,i,e}=0$ if served elsewhere), and arrival $j$ incurs execution requirement $C_{t,j,e}$ if admitted to endpoint $e$. We enforce memory capacity at epoch boundaries as a discrete-time approximation: an initially active turn and a newly admitted turn occupy boundary memory $K_{t,i,e}^{\mathrm{act}}$ and $K_{t,j,e}^{\mathrm{adm}}$ at endpoint $e$, respectively (equal to zero if the turn completes during the epoch), whereas completed turns consume boundary memory only through their retained idle caches of observed size $b_i(t)$.

All per-session execution requirements and active-memory footprints are nonnegative and uniformly bounded. Conditional on the beginning-of-epoch state and arrival sequence, these random quantities are independent across sessions and follow known distributions determined by each session's state and action, though their realizations are unknown at the time of routing.

A nonanticipative policy $\pi$ specifies the routing actions $e_{t,j}$ and retained cache sets $\mathcal I_e(t)$ using only the information observed up to each decision time. Letting $N^{\mathrm{cmp}}(t)$ denote the number of turns completed during epoch $t$, our objective from initial state $M(0)=M^0$ is to maximize expected completions:
\begin{equation}
\label{eq:stochastic-control-problem}
\begin{aligned}
  V_0(M^0)=\sup_{\pi}\quad&
  \E^\pi\!\left[\sum_{t=0}^{H-1}N^{\mathrm{cmp}}(t)\right]\\
  \text{subject to}\quad&
  \sum_{i\in\mathcal A_t}L_{t,i,e}
  +\sum_{j=1}^{J_t}C_{t,j,e}\mathbf 1\{e_{t,j}=e\}
  \le B_e,\\
  &\sum_{i\in\mathcal A_t}K_{t,i,e}^{\mathrm{act}}
  +\sum_{j=1}^{J_t}K_{t,j,e}^{\mathrm{adm}}
     \mathbf 1\{e_{t,j}=e\}
  +\sum_{i\in\mathcal I_e(t)}b_i(t)
  \le G_e,
\end{aligned}
\end{equation}
where both constraints must hold almost surely for every $e\in\cE$ and $t=0,\ldots,H-1$. The first constraint bounds computation within epoch $t$, while the second allocates boundary memory between unfinished turns and retained idle caches.

The objective also admits a useful accounting identity. Define $N_A(M)=\sum_{s\in\cSA}M_s$ as the total active population. Because every admitted turn either completes during the epoch or remains active at the next boundary,
\[
  N^{\mathrm{cmp}}(t)
  =N_A(M(t))
   +\sum_{j=1}^{J_t}\mathbf 1\{e_{t,j}\ne0\}
   -N_A(M(t+1)).
\]
Summing over $t=0,\ldots,H-1$ yields
\[
  V_0(M^0)=N_A(M^0)+\sup_\pi\E^\pi\!\left[
    \sum_{t=0}^{H-1}\sum_{j=1}^{J_t}\mathbf 1\{e_{t,j}\ne0\}
    -N_A(M(H))\right].
\]
Maximizing completions is therefore equivalent to maximizing admissions minus the unfinished turns remaining at the terminal time $H$.

Solving \eqref{eq:stochastic-control-problem} requires trading off immediate service against future cache reuse: routing determines where subsequent turns can reuse computations, while retention determines whether those caches remain available. Because these decisions interact across sessions through shared computation and memory limits, exact dynamic programming is intractable. The next section develops a fluid approximation that preserves these intertemporal interactions and leads to a tractable policy.

\begin{remark}[Emergency intervention]\label{remark:emergency}
The formulation in \eqref{eq:stochastic-control-problem} describes normal operation. For theoretical analysis, we assume that a causal emergency rule can halt admissions and interrupt active turns before any endpoint exceeds capacity, with interrupted turns earning no completion reward. This rule serves solely to keep the policy well defined on exceptional sample paths and is not intended for practical implementation; our policy reserves capacity to make such interventions rare, and we bound their impact on performance separately. 
\end{remark}

\section{The Fluid Approximation and Structural Insights}\label{sec:fluid}
{
We develop a fluid approximation to characterize how anticipated cache reuse should govern current admission, routing, and retention decisions. Routing a turn determines where a conversation's history can next be reused, whereas retention determines whether that history survives until the conversation returns. Theorem~\ref{prop:bp-control} couples these decisions through endpoint-specific shadow prices for computation and memory together with each conversation's continuation value.

For routing, the fluid model favors the endpoint with the highest bid, combining the value of serving the current turn, the cost of the resources it consumes, and the resulting change in continuation value. When endpoints face similar congestion, keeping a conversation at its cached endpoint avoids rebuilding its prefix. By contrast, a persistent congestion imbalance can justify migrating a conversation even when doing so makes the current turn more expensive: the rebuild penalty is incurred once, whereas lower congestion yields recurring savings over subsequent turns. Consequently, two conversations facing the same load imbalance may warrant different routes depending on their remaining turn horizons.

For retention, the fluid model preserves an idle cache whenever the continuation value it protects exceeds the one-epoch opportunity cost of the memory it occupies. This trade-off explains why intermediate-length prefixes can deserve priority: short prefixes are inexpensive to rebuild, whereas very long prefixes may consume more memory than a conversation's remaining turns can justify. Return timing matters as well: when tool execution takes a predictable amount of time, a session that has already been idle for several epochs may be the closest to returning.

These structural insights directly guide the design of the practical algorithm \textsc{PARK} in Section~\ref{sec:swe-gradient-implementation} later. Its admission and routing controller updates resource prices online and approximates the continuation benefit of migration by comparing the cost of one subsequent warm turn across endpoints. Meanwhile, its cache manager ranks idle prefixes by their estimated reuse benefit per unit of memory and time, accounting for both prefix length and the effective waiting time per successful reuse.
}

\subsection{The Fluid Linear Program}
Replacing random quantities with their expectations and discrete
decisions with continuous fractions does not immediately yield a
linear program. To see why, let $m_s(t)$ denote the fluid population
in state $s$ (with $m_s(0)=M_s^0$), and let $\alpha_{s,e}(t)$ be
the fraction of arriving state-$s$ prompts routed to endpoint $e$.
The resulting execution load is
\(
  c_{s,e}(t)\bigl[p_s(t)m_s(t)+\nu_s(t)\bigr]\alpha_{s,e}(t),
\)
where $\nu_s(t)$ is the expected number of new state-$s$ arrivals
and $c_{s,e}(t)$ is the conditional mean execution requirement of a
state-$s$ prompt admitted to $e$. Because both the future
population $m_s(t)$ and the routing fraction $\alpha_{s,e}(t)$
depend on policy decisions, their product is bilinear. An
analogous nonlinearity arises when multiplying the post-service
idle population by the fraction of caches retained.

To restore linearity, we instead optimize directly over state--action
flow volumes: the total mass of prompts routed and caches retained
in each state. This reformulation makes both resource consumption
and population dynamics linear while preserving how current
decisions govern future cache availability. We first present the
resulting linear program, and then use its dual prices to derive
routing and retention principles that account for a conversation's
future turns.

\paragraph{Arrival and routing flows.}
Let $y_{s,e}(t)$ denote the volume of state-$s$ prompts assigned to
action $e\in\cEzero$ in epoch $t$, where $e=0$ represents rejection.
Arrival conservation requires
\begin{equation}
  \sum_{e\in\cEzero}y_{s,e}(t)
  =p_s(t)m_s(t)+\nu_s(t),
  \qquad s\in\cSI.
  \label{eq:fluid-routing-conservation}
\end{equation}
Because new sessions enter without a cache, $\nu_s(t)=0$ whenever
$o(s)\ne0$. The remaining population $(1-p_s(t))m_s(t)$ consists
of idle sessions that submit no prompt during epoch~$t$.

\paragraph{Population after service.}
Let $\widetilde m_{s'}(t+1)$ be the pre-retention population in live
state $s'$ at the end of epoch~$t$, after service and session
transitions. Applying the transition kernels from Section~\ref{sec:model} gives
\begin{align}
  \widetilde m_{s'}(t+1)
  ={}&
  \sum_{s\in\cSA}m_s(t)Q_t^A(s,s')
  +\sum_{s\in\cSI}(1-p_s(t))m_s(t)Q_t^I(s,s')
  \nonumber\\
  &+\sum_{s\in\cSI}\sum_{e\in\cEzero}
       y_{s,e}(t)P_t^e(s,s'),
  \qquad s'\in\cS.
  \label{eq:fluid-pre-retention}
\end{align}
The three sums track transitions from initially active sessions,
unprompted idle sessions, and sessions with arriving prompts,
respectively.

\paragraph{Cache retention and the next population.}
For each idle state $s\in\cSI$, let $z_{s,\KEEP}(t)$ and
$z_{s,\DROP}(t)$ denote the flow volumes retained and evicted at
the epoch boundary, satisfying
\begin{equation}
  z_{s,\KEEP}(t)+z_{s,\DROP}(t)
  =\widetilde m_s(t+1),
  \qquad s\in\cSI.
  \label{eq:fluid-retention-conservation}
\end{equation}
Cache-free states satisfy $z_{s,\KEEP}(t)=0$, so all mass is
assigned to eviction and remains unchanged. Recalling that $\Ev(s)$
maps an idle state $s$ to its cache-free counterpart with $o=0$, the
post-retention population is
\begin{equation}
  m_{s'}(t+1)=
  \begin{cases}
    \widetilde m_{s'}(t+1),
      &s'\in\cSA,\\[3pt]
    z_{s',\KEEP}(t)
    +\displaystyle\sum_{\substack{s\in\cSI:\\\Ev(s)=s'}}
       z_{s,\DROP}(t),
      &s'\in\cSI.
  \end{cases}
  \label{eq:fluid-post-retention}
\end{equation}
Retention thus preserves a session's cache location, whereas
eviction shifts its mass into the corresponding cache-free state.

\paragraph{Resource constraints.}
Let $\ell_{s,e}(t)$ denote the conditional mean execution
requirement of an initially active state-$s$ session at endpoint $e$,
and recall that $c_{s,e}(t)$ is the conditional mean execution
requirement of a state-$s$ prompt admitted to $e$. Per-epoch
execution capacity requires
\begin{equation}
  \sum_{s\in\cSA}\ell_{s,e}(t)m_s(t)
  +\sum_{s\in\cSI}c_{s,e}(t)y_{s,e}(t)
  \le B_e,
  \qquad e\in\cE.
  \label{eq:fluid-execution}
\end{equation}

Similarly, let $\kappa_{s,e}^{\mathrm{act}}(t)$ and
$\kappa_{s,e}^{\mathrm{adm}}(t)$ denote the conditional mean
boundary-memory footprints of an initially active state-$s$ turn at
endpoint $e$ and a state-$s$ prompt admitted to $e$, respectively
(both equal zero if the turn completes during the epoch). Letting
$b_s$ denote the idle cache size recorded in state $s$, boundary
memory capacity requires
\begin{equation}
  \sum_{s\in\cSA}\kappa_{s,e}^{\mathrm{act}}(t)m_s(t)
  +\sum_{s\in\cSI}\kappa_{s,e}^{\mathrm{adm}}(t)y_{s,e}(t)
  +\sum_{\substack{s\in\cSI:\\o(s)=e}}b_s z_{s,\KEEP}(t)
  \le G_e,
  \qquad e\in\cE.
  \label{eq:fluid-memory}
\end{equation}
The three terms account for initially active turns, newly admitted
turns, and retained idle caches, mirroring the stochastic memory
constraint in \eqref{eq:stochastic-control-problem}.

\paragraph{The remaining-horizon optimization problem.}
Consider the remaining-horizon planning problem starting at epoch $t$
with population vector $M$. Applying the completion accounting
identity from Section~\ref{sec:model}, we maximize total admissions minus the
terminal active population:
\begin{equation}
\label{eq:fluid-lp-t}
\begin{aligned}
  \mathcal V_t(M)
  =\max_{m,\widetilde m,y,z}\quad&
  \sum_{\tau=t}^{H-1}\sum_{s\in\cSI}\sum_{e\in\cE}
      y_{s,e}(\tau)-N_A(m(H))\\
  \text{subject to}\quad
  &m(t)=M,\\
  &\eqref{eq:fluid-routing-conservation}
      \text{--}\eqref{eq:fluid-memory},
      &&\tau=t,\ldots,H-1,\\
  &z_{s,\KEEP}(\tau)=0,
      &&s\in\cSI:\ o(s)=0,\quad \tau=t,\ldots,H-1,\\
  &m,\widetilde m,y,z\ge0.
\end{aligned}
\end{equation}
Total fluid completions equal $N_A(M)+\mathcal V_t(M)$, and setting
$t=0$ and $M=M^0$ recovers the initial planning problem.

This linear program links current routing and retention decisions
to future service requirements through the evolving distribution of
cache locations. Its dual shadow prices allow us to weigh the
immediate resource cost of each action against the continuation
value of the session state it creates, and we formalize this comparison
next to derive routing and retention principles.

\subsection{Bid Prices and the Value of Future Cache Reuse}
\label{sec:bp-insights}
Fix a feasible remaining-horizon LP \eqref{eq:fluid-lp-t} starting at epoch $t_0$ with population $M$, and let
\[
  \lambda_e(t),\gamma_e(t)\ge0,
  \qquad e\in\cE,\quad t=t_0,\ldots,H-1,
\]
denote its optimal dual shadow prices for execution and memory capacity, respectively. All prices and values below refer to this same LP and may vary across epochs.

At these prices, let $v_s(t+1)$ denote the continuation value of a session in state $s$ at the start of epoch $t+1$, net of future dual resource charges and inclusive of the terminal penalty for unfinished turns in \eqref{eq:fluid-lp-t}. Let $w_s(t)$ denote its pre-retention value at the end of epoch $t$; Appendix~\ref{app:bp-dual} gives the backward recursion defining both value functions. For an idle state $s$ with cache owner $e$, let $s^0=\Ev(s)$ be the corresponding cache-free state. Retaining the cache yields continuation gain $\Delta_s(t+1)=v_s(t+1)-v_{s^0}(t+1)$ at boundary memory cost $\gamma_e(t)b_s$, whereas eviction transitions the session to $s^0$. Consequently,
\begin{equation}
  w_s(t)=\max\bigl\{v_{s^0}(t+1),\,
                   v_s(t+1)-\gamma_e(t)b_s\bigr\}.
  \label{eq:bp-boundary-option}
\end{equation}
For active and cache-free idle states, $w_s(t)=v_s(t+1)$.

\paragraph{The bid-price characterization.}
For a state-$s$ prompt routed to endpoint $e\in\cE$, define its immediate dual resource cost by
\begin{equation}
  r_{s,e}(t)=\lambda_e(t)c_{s,e}(t)
    +\gamma_e(t)\kappa_{s,e}^{\mathrm{adm}}(t).
  \label{eq:bp-resource-rent}
\end{equation}
The net value of admitting the prompt to endpoint $e$ relative to rejection, termed its \emph{bid}, is
\begin{equation}
  \mathcal B_{s,e}(t)
  =1-r_{s,e}(t)
   +\sum_{s'\in\cS}
       \bigl[P_t^e(s,s')-P_t^0(s,s')\bigr]w_{s'}(t),
  \qquad e\in\cE,
  \label{eq:bp-bid}
\end{equation}
with $\mathcal B_{s,0}(t)=0$. The three terms in \eqref{eq:bp-bid} capture the unit admission reward, the immediate dual resource cost, and the net change in expected continuation value relative to rejection.

The next theorem characterizes optimal fluid routing and retention decisions in terms of these bids and shadow prices; its proof appears in Appendix~\ref{app:bp-control}.
\begin{theorem}[Optimal routing and retention rules]
\label{prop:bp-control}
For any optimal fluid solution and the dual prices above:
\begin{enumerate}[label=\alph*.]
  \item A positive routing flow $y_{s,e}(t)>0$ uses an action maximizing $\mathcal B_{s,e}(t)$ over $\cEzero$.
  \item For a cached idle state with positive pre-retention population, its cache is retained if $\Delta_s(t+1)>\gamma_{o(s)}(t)b_s$ and evicted if the inequality is reversed.
\end{enumerate}
At equality, the flow and capacity constraints determine the feasible allocation among tied actions.
\end{theorem}

\paragraph{When routing should follow the existing cache.}
For any two endpoints $e_1,e_2\in\cE$, define their continuation-value differential by
\(
  \Theta_{s;e_1,e_2}(t)
  =\sum_{s'\in\cS}
    \bigl[P_t^{e_1}(s,s')-P_t^{e_2}(s,s')\bigr]w_{s'}(t).
\)
Comparing bids via \eqref{eq:bp-bid} gives
\begin{equation}
  \mathcal B_{s,e_1}(t)\ge\mathcal B_{s,e_2}(t)
  \quad\Longleftrightarrow\quad
  \Theta_{s;e_1,e_2}(t)\ge r_{s,e_1}(t)-r_{s,e_2}(t).
  \label{eq:bp-owner-threshold}
\end{equation}
Thus, routing to $e_1$ rather than $e_2$ is optimal whenever its continuation-value advantage outweighs any higher immediate resource cost.

Suppose endpoint $e_2$ holds the session's cache. When shadow prices and future service conditions are similar across endpoints, cache reuse lowers $c_{s,e_2}(t)$ and makes $e_2$ cheaper to use, favoring warm routing. However, if $e_2$ becomes sufficiently congested, its higher shadow prices raise $r_{s,e_2}(t)$ and depress future continuation values there, making migration to a less loaded endpoint $e_1$ worthwhile despite the upfront cost of rebuilding the cache. This trade-off mirrors the intuition behind SGLang's \texttt{cache\_aware} routing rule~\cite{sglang-cache-aware}, which prioritizes cache locality when server loads are balanced but switches to the least-loaded endpoint once load imbalance exceeds a threshold.

Condition~\eqref{eq:bp-owner-threshold} also highlights the role of remaining conversation length. Serving the current turn at $e_1$ establishes a warm cache there for subsequent turns; if $e_1$ is expected to remain less congested, those future turns yield recurring savings that amortize the one-time rebuild cost. Migration is therefore attractive for a conversation with many turns remaining, whereas it will not pay off for a session nearing completion. Consequently, the congestion gap needed to justify rerouting depends jointly on expected remaining turns, the persistence of the destination's cost advantage, and the memory price of retaining the newly built cache.

While routing determines where the next cache is created, retention determines whether it remains available for reuse. The next subsection uses Theorem~\ref{prop:bp-control}.b to examine how cache length, idle age, and return timing govern which caches should be retained.

\subsection{Which Idle Caches Should Be Retained?}
\label{sec:bp-retention}
An idle cache occupies memory now to save computation when the next
prompt arrives. Theorem~\ref{prop:bp-control}.b compares this
future benefit against the boundary memory charge: retain the cache
when $\Delta_s(t+1)>\gamma_{o(s)}(t)b_s$, and evict it when the
inequality is reversed. We use this comparison to examine how
\emph{cache length}, \emph{return timing}, and \emph{idle age}
govern the value of retaining an idle cache.

Fix an idle cache at endpoint $e$ at the boundary ending epoch $t$,
with prefix length $L>0$ and memory footprint $b_s=\beta L$
($\beta>0$). There are $N=H-t-1\ge1$ remaining prompt epochs.
Let $T$ denote the number of epochs until the next prompt, with
$T=\infty$ if the session does not return within the horizon. For
tractability, suppose the endpoint's memory shadow price remains
constant at $\gamma>0$ across waiting boundaries;
Appendix~\ref{app:bp-retention} provides exact conditions and
extends the analysis to time-varying prices.

\paragraph{The cost of waiting for reuse.}
Suppose we retain the cache for at most $u$ epochs, evicting it
if no prompt arrives by then. Over this waiting window, the expected
storage duration is $\E[\min(T,u)]$ and the probability of reusing
the cache before eviction is $\Pp(T\le u)$. Define the
\emph{effective waiting time} per cached return by
\begin{equation}
 \deff=
 \min_{\substack{1\le u\le N\\\Pp(T\le u)>0}}
 \frac{\E[\min(T,u)]}{\Pp(T\le u)},
 \label{eq:bp-effective-return}
\end{equation}
with $\deff=\infty$ when no return is possible (in which case
immediate eviction is optimal). Thus, $\deff$ measures the minimum
expected storage time per successful cache reuse across all
stopping deadlines $u$. For example, if a prompt arrives
deterministically in five epochs, then $\deff=5$, so retaining the
cache can be optimal even when an immediate return in the next
epoch has probability zero. For random return times, $\deff$
incorporates the option to abandon waiting early and therefore
differs from the unconditional mean $\E[T]$.

\paragraph{Conditions for keeping or evicting a cache.}
For a session with a positive probability of return, let $R<\infty$
uniformly bound its expected remaining completions under any
admissible policy, conditional on the history and chosen route at a
cached return or at the first cold admission after eviction.
Suppose further that, following eviction, the first cold admission
incurs an expected execution charge of at least
\begin{equation}
 \underline\lambda(a_1L+a_2L^2),
 \qquad a_1\ge0,\quad a_2>0,\quad\underline\lambda\ge0,
 \label{eq:bp-retention-cold-bound}
\end{equation}
where $\underline\lambda$ lower-bounds the execution shadow price
and the linear and quadratic terms capture the computation required
to rebuild the lost prefix (see, e.g.,
\citealt{vaswani,hoffmann-compute}). When the session instead
returns with its cache intact, its continuation value equals
expected future completions minus future execution and memory
charges; let $\varepsilon(L)$ uniformly bound how far this value
falls below $R$ across all feasible return dates. A smaller
$\varepsilon(L)$ thus indicates that retaining the cache preserves
more of the session's potential completion value
(Appendix~\ref{app:bp-retention-deficit} formalizes
$\varepsilon(L)$ and provides explicit upper bounds).

The following proposition establishes sufficient conditions for
fluid retention and eviction, showing how cache length and return
timing jointly govern eviction; its proof appears in
Appendix~\ref{app:bp-retention-proof}.
\begin{proposition}
\label{prop:bp-retention-length}
Under the preceding setting, the following are sufficient conditions
for the fluid retention decision.
\begin{itemize}
    \item For long cache,
    \begin{equation}
        \gamma b_s\deff
 >R
 \quad\Longrightarrow\quad\text{evict}.
 \label{eq:bp-retention-long-drop}
    \end{equation}
    \item For cache with intermediate length,
    \begin{equation}
         \gamma b_s\deff+\varepsilon(L)
 <\min\{\underline\lambda(a_1L+a_2L^2),R\}
 \quad\Longrightarrow\quad\text{retain}.
 \label{eq:bp-retention-keep}
    \end{equation}

    \item For short cache, additionally suppose that,
under each fixed routing action, cache reuse changes only
execution cost, with savings at most
$\bar\lambda_e(a_1L+a_2L^2)$, where $\bar\lambda_e\ge0$.
Then
\begin{equation}
 \bar\lambda_e(a_1+a_2L)<\beta\gamma\deff
 \quad\Longrightarrow\quad\text{evict}.
 \label{eq:bp-retention-short-drop}
\end{equation}
\end{itemize}
\end{proposition}

Condition~\eqref{eq:bp-retention-keep} compares the minimum storage
cost per cached return against a lower bound on the net value
preserved by the cache, whereas \eqref{eq:bp-retention-long-drop}
and \eqref{eq:bp-retention-short-drop} reverse this comparison
using upper bounds on that benefit.
Proposition~\ref{prop:bp-retention-length} yields three structural
implications.

\paragraph{Why intermediate lengths can be worth protecting.}
\begin{itemize}
 \item \textbf{Short caches save too little computation.}
 When prefix-recomputation savings are primarily quadratic in
 length ($a_1\approx 0$), reuse savings shrink faster than linear
 storage costs as $L$ decreases. Under the additional service
 assumption, \eqref{eq:bp-retention-short-drop} implies that
 sufficiently short caches are cheaper to rebuild than to retain,
 holding return prospects and shadow prices fixed.
 \item \textbf{Very long caches consume too much memory.}
 By \eqref{eq:bp-retention-long-drop}, retaining a cache is
 suboptimal whenever $L>R/(\beta\gamma\deff)$. Even a costly rebuild
 cannot justify storage if the remaining conversation offers too few
 expected completions $R$ to amortize the memory rent.
 \item \textbf{Intermediate-length caches merit highest retention priority.}
 Caches of moderate length are expensive enough to rebuild yet
 compact enough to store profitably.
 Appendix~\ref{app:bp-retention-primitive} derives an explicit
 interval of protected cache lengths under simple service-cost
 bounds.
\end{itemize}

\paragraph{How waiting for reuse affects retention.}
Cache value also depends on how long memory must be held before
the next prompt arrives. Holding cache length, continuation bounds,
and shadow prices fixed, a larger $\deff$ increases the storage cost
per return, tightening the retention condition
\eqref{eq:bp-retention-keep} and relaxing the eviction conditions
\eqref{eq:bp-retention-long-drop}--\eqref{eq:bp-retention-short-drop}.
In particular, the maximum justifiable cache length
$R/(\beta\gamma\deff)$ shrinks as $\deff$ grows, so even a cache
that is expensive to rebuild should be evicted if capturing its
reuse requires excessive waiting.

\paragraph{When idle age is a useful guide.}
An older cache is not automatically less valuable; what matters is
what its idle age reveals about the residual time until the next
prompt. For a given session class, let $U$ denote the total idle
duration (in epochs) from the previous turn's completion to the next
prompt, possibly infinite. For a session already idle for $a$
epochs, define its conditional return probability (hazard rate) in
the next epoch by
\[
 p(a,L)=\Pp(U=a+1\mid U>a,L).
\]
The following corollary characterizes how the effective waiting time $\deff(a,L)$ varies with a session's idle age $a$ under decreasing, increasing, and constant return-probability regimes.
\begin{corollary}
\label{cor:bp-retention-age}
Compare idle ages at a fixed cache length and remaining horizon.
If $p(a,L)$ is nonincreasing in $a$, then
\[
 \deff(a,L)=1/p(a,L),
\]
with $1/0=\infty$, and $\deff$ is nondecreasing in age.
If $p(a,L)$ is nondecreasing in $a$, then $\deff$ is
nonincreasing in age. If $p(a,L)$ is constant in $a$, then
$\deff$ is independent of age.
\end{corollary}

Because a larger $\deff$ favors eviction,
Corollary~\ref{cor:bp-retention-age} explains when
\emph{least-recently-used} (LRU; see, e.g., \citealt{mattson1970storage,zheng2024sglang}) eviction is
effective and when it fails. When the return hazard $p(a,L)$
decreases with idle age, meaning longer inactivity signals a less
imminent return, $\deff$ increases with $a$, making LRU sensible.
Conversely, when a prompt is expected after a bounded processing or
tool-execution delay, the return hazard increases with age: an older
session is closer to returning, so LRU can prematurely evict a
cache that is about to be reused. Combining cache length with
return-time behavior therefore provides a sharper eviction
criterion than recency alone.

The routing and retention principles developed in this section
characterize optimal decisions in the deterministic fluid model.
Implementing these insights in a stochastic system requires
safeguarding against random fluctuations in arrivals, service
durations, and memory footprints. The next section introduces a
buffer-reserve policy that translates the fluid solution into
online decisions and proves that its performance loss relative to
the fluid benchmark grows sublinearly with system scale.

\section{A Fluid-Based Policy: Theory and Implementation}
\label{sec:br-policy}
{In this section, we first present a buffer-reserve policy that implements the fluid solution and establish its asymptotic performance guarantee, and we then adapt its decision rules for practical cluster deployment. The buffer-reserve policy inherits the fluid model's bid-price structure. For routing, an endpoint's bid weighs the value of serving a turn, including future cache reuse, against its shadow cost of computation and memory. For retention, an idle cache is kept whenever its continuation benefit exceeds the shadow cost of the memory it occupies. Building on this logic, our practical algorithm, \textsc{PARK} (described in Section~\ref{sec:swe-gradient-implementation}), approximates these comparisons using online resource prices and a one-turn lookahead for admission and routing, alongside an index of expected reuse benefit per unit of memory and time for cache retention.}

\subsection{The Buffer-Reserve Policy}

The fluid LP plans resource allocation using conditional expectations,
whereas actual arrivals, execution requirements, and memory
footprints fluctuate around their means. Allocating full nominal
capacity can therefore trigger an overload within the current epoch
that re-solving at the next epoch cannot undo. Multi-turn service
compounds this risk across epochs: unfinished turns carry inherited
resource commitments into future epochs that cannot be shed by
rejecting new prompts, while random completions govern when and
where idle caches appear. To guard against both within-epoch
spikes and multi-turn population drift, we pair progressive capacity
reserves with epoch-by-epoch reoptimization.

Fix a reserve parameter $r>0$. When solving the fluid LP at epoch
$t$, we reduce the available capacities at each remaining epoch
$\tau$ to
\begin{equation}
 B_e^{-,t}(\tau)=B_e-(\tau-t+1)r,
 \qquad
 G_e^{-,t}(\tau)=G_e-(\tau-t+1)r,
 \qquad \tau=t,\ldots,H-1.
 \label{eq:br-reserved-capacities}
\end{equation}
Thus, epoch $t$ holds back $r$ units of each resource, and the
buffer grows linearly across future epochs. When the LP is re-solved
at $t+1$ from the observed population $M(t+1)$, each future epoch
$\tau\ge t+1$ regains $r$ units of capacity relative to the
previous plan. This released slack ensures recourse feasibility
when realized populations, particularly ongoing active turns, exceed
their fluid predictions, while reoptimization adapts routing and
retention quotas to realized cache locations.

Algorithm~\ref{def:br-policy} converts the current-epoch fluid
flows into rounded integer quotas, making immediate routing
decisions upon arrival and boundary retention decisions after
service outcomes are observed. Whenever a capacity intervention is
required or the re-solved LP becomes infeasible, the causal
emergency rule from Section~3 applies.

\begin{algorithm}[t]
\caption{Buffer-reserve fluid policy}
\label{def:br-policy}
At each epoch $t$, until an emergency occurs:
\begin{enumerate}[leftmargin=*,itemsep=3pt,topsep=3pt]
\item \textbf{Plan.}
Solve \eqref{eq:fluid-lp-t} from $M(t)$ with capacities
\eqref{eq:br-reserved-capacities}. If infeasible, invoke the
emergency rule. Otherwise, select an optimal solution via a fixed
deterministic tie-breaking rule, and denote its current-epoch flows
by $y_{s,e}(t)$ and $z_{s,\KEEP}(t)$.

\item \textbf{Route.}
Upon arrival $j$ in state $s_{t,j}$, assign $e_{t,j}=e$ to the
first endpoint $e\in\cE$ in a fixed state-specific order that has
remaining quota:
\[
 \sum_{i=1}^{j-1}
 \1\{s_{t,i}=s_{t,j},\ e_{t,i}=e\}
 <\lfloor y_{s_{t,j},e}(t)\rfloor.
\]
If no such endpoint exists, reject the prompt ($e_{t,j}=0$).
Unused quotas expire at the epoch boundary.

\item \textbf{Retain.}
Observe the pre-retention states $\widetilde s_i(t+1)$, including
newly completed sessions. For each idle state $s$ with
$o(s)=e\ne0$, retain $\lfloor z_{s,\KEEP}(t)\rfloor$ caches (or
all available caches if fewer exist) in a fixed session order,
forming $\mathcal I_e(t)$. Evict all remaining idle caches, update
$M(t+1)$, and proceed to epoch $t+1$.
\end{enumerate}
\end{algorithm}

\paragraph{Connection to bid-price control.}
By Theorem~\ref{prop:bp-control}, positive flows in the buffered LP
follow the optimal bid-price routing and retention comparisons
under the buffered dual prices, while the flow magnitudes resolve
ties among indifferent actions and enforce capacity limits.
Algorithm~\ref{def:br-policy} implements these primal--dual
decisions through rounded quotas that are refreshed at every
reoptimization.

\subsection{Scaling and Performance Guarantee}
\label{sec:br-scaling}

To evaluate the policy in large-scale deployments where servers house multiple GPUs and support many concurrent sessions, we analyze an asymptotic regime in which demand and capacity grow proportionally with a scale parameter $n\ge2$. Specifically, in system~$n$, endpoint capacities and expected new-session arrivals scale as
\begin{equation}
 B_e^{(n)}=nB_e,\qquad
 G_e^{(n)}=nG_e,\qquad
 \E[\Xi_s^{(n)}(\tau)\mid\mathcal F_\tau]
 =n\nu_s(\tau),
 \label{eq:br-scaling}
\end{equation}
where $\Xi_s^{(n)}(\tau)$ counts new state-$s$ sessions entering in epoch~$\tau$ and $\mathcal F_\tau$ denotes the history observed up to the start of epoch~$\tau$. The set of endpoints $\cE$, state space $\cS$, horizon $H$, and all per-session resource and transition laws remain fixed, while normalized initial populations lie in a bounded set $\|m\|_1\le K_0$ independent of~$n$.

The following assumption bounds tail fluctuations in new-session arrivals.
\begin{assumption}
\label{ass:br-scaling}
For every state $s$, epoch $\tau$, and $u\ge0$,
\[
 \Pp\!\left(
 |\Xi_s^{(n)}(\tau)-n\nu_s(\tau)|>u
 \,\middle|\,\mathcal F_\tau
 \right)
 \le C_0\exp\!\left\{-\frac{u^2}{C_0n}\right\},
\]
where $C_0>0$ is independent of $n$ and the observed history.
\end{assumption}
This sub-Gaussian concentration holds, for instance, when arrivals originate from $O(n)$ independent user sources, each generating a bounded number of new sessions per epoch conditional on $\mathcal F_\tau$, in which case Hoeffding's inequality yields the bound.

In system~$n$, we set the capacity reserve to
\begin{equation}
 r_n=b\sqrt{n\log n},
 \label{eq:br-buffer-size}
\end{equation}
for a fixed constant $b>0$, and apply Algorithm~\ref{def:br-policy} with scaled capacities $nB_e$ and $nG_e$ and reserve $r=r_n$ in \eqref{eq:br-reserved-capacities}. Starting at epoch~$t$ from population $M=nm$, let $V_{t,n}^{\mathrm{BR}}(nm)$ denote the policy's expected completions through horizon~$H$, and let $\mathsf E_n$ be the event that an emergency intervention occurs at any point over this period (including the terminal boundary). We benchmark performance against the unbuffered fluid value
\[
 \mathcal W_t(m)=N_A(m)+\mathcal V_t(m),
\]
where $\mathcal V_t(m)$ is the optimal value of \eqref{eq:fluid-lp-t} under nominal capacities $B_e,G_e$ and arrival rates $\nu_s(\tau)$ without reserves. By linearity, $n\mathcal W_t(m)$ is the fluid benchmark for system~$n$. The following theorem shows that the buffer-reserve policy avoids capacity emergencies with high probability and tracks this fluid benchmark within $\widetilde O(\sqrt n)$ completions.

\begin{theorem}
\label{thm:br-regret}
Under Assumption~\ref{ass:br-scaling}, choose $b$ sufficiently large. Then for all sufficiently large $n$ and any normalized initial population $m$ for which the initial buffered LP is feasible:
\begin{enumerate}[label=(\roman*)]
 \item The policy is nonanticipative, satisfies all resource capacities (with the emergency rule), and incurs an emergency intervention with probability at most $\Pp^{\mathrm{BR}}(\mathsf E_n)\le n^{-2}$.

 \item Its expected completions satisfy
 \begin{equation}
   n\mathcal W_t(m)-V_{t,n}^{\mathrm{BR}}(nm)
   \le C\sqrt{n\log n},
   \label{eq:br-regret}
 \end{equation}
 where $C>0$ is independent of $n$ and $m$.
\end{enumerate}
\end{theorem}

Theorem~\ref{thm:br-regret} establishes an $\widetilde O(\sqrt n)$ additive shortfall relative to the fluid benchmark; dividing \eqref{eq:br-regret} by $n$ gives
\[
 \mathcal W_t(m)-\frac{V_{t,n}^{\mathrm{BR}}(nm)}{n}
 \le
 C\sqrt{\frac{\log n}{n}},
\]
so the per-unit performance loss vanishes as $n$ grows. 
The proof in Appendix~\ref{app:br-proofs} decomposes the performance gap into three components: quota rounding and stochastic fluctuations in arrivals and transitions contribute $O(\sqrt n)$, the capacity reserve $r_n$ sacrifices $O(\sqrt{n\log n})$, and rare emergency interventions contribute $O(1/n)$. Setting $r_n=b\sqrt{n\log n}$ balances these effects: it provides enough slack at each reoptimization to absorb drift in unfinished work and preserve LP feasibility with high probability, while adapting routing and retention quotas to realized cache locations.

\begin{remark}[{Asymptotic optimality}]
\label{rem:br-optimality}
{
Theorem~\ref{thm:br-regret} can also yield an asymptotic optimality
guarantee when all policies have access to the same emergency rule.
For this comparison, we charge each policy a fixed penalty per
interrupted turn, large enough to cover the benefit of freeing
capacity. The emergency rule is a theoretical safeguard: the
buffer-reserve policy reserves enough capacity that it needs to
interrupt service only on rare sample paths. Its expected total
penalty therefore vanishes as the system grows, and its expected
completions net of penalties are within $O(\sqrt{n\log n})$ of the
best feasible policy. Appendix~\ref{app:br-proofs} gives the proof.
}
\end{remark}

{
\subsection{Practical Implementation: Asynchronous Primal--Dual Control and Index Eviction}
\label{sec:swe-gradient-implementation}

The buffer-reserve policy (Algorithm~\ref{def:br-policy}) operates in synchronous epochs, observes the full system state before each decision, re-solves a multi-epoch linear program at every step, and directly enforces both routing and cache eviction. While these abstractions facilitate theoretical analysis, production LLM serving clusters differ in three operational respects:
\begin{enumerate}[label=(\roman*)]
  \item \emph{Asynchronous arrivals and delayed state updates:} Prompts arrive continuously between control updates, while endpoint utilization snapshots arrive with delay, making it computationally prohibitive to re-solve a multi-epoch LP at every arrival.
  \item \emph{Decoupled admission and dispatch:} Because routing decisions precede actual client dispatch, the controller must reserve capacity through short-lived \emph{admission offers} to prevent concurrent requests from double-booking resources.
  \item \emph{Engine-local cache reclamation:} In production inference engines (such as vLLM; see, e.g., \citealt{vllm-prefix}), the central router does not forcibly evict KV blocks; instead, each endpoint's local memory manager reclaims unpinned idle blocks on demand whenever active generation requires space.
\end{enumerate}

To bridge these gaps while preserving the economic structure of Theorem~\ref{prop:bp-control} and Algorithm~\ref{def:br-policy}, we pair an event-driven admission and routing controller (approximating the bid-price rule in Theorem~\ref{prop:bp-control}.a) with an endpoint-local retention index (approximating the retention rule in Theorem~\ref{prop:bp-control}.b). Rather than re-solving the full fluid LP, the central controller maintains shadow prices $\lambda_e \ge 0$ for computation and $\gamma_e \ge 0$ for active memory at each endpoint $e \in \cE$, uses them to evaluate bid scores, and updates them via projected subgradient steps on excess desired demand. Simultaneously, each endpoint ranks reclaimable idle caches using a Gittins-style value-per-byte-epoch index derived from Section~\ref{sec:bp-retention}. We refer to the resulting integrated policy as \textsc{PARK} (summarized in Algorithm~\ref{alg:online-dual}), with experimental parameter values listed in Table~\ref{tab:online-dual-parameters}, and full implementation details included in Section~\ref{app:park-implementation}, both in the Appendix. Table~\ref{tab:br-park} contrasts how \textsc{PARK} and the buffer-reserve policy translate the two bid-price comparisons of Theorem~\ref{prop:bp-control} into serving decisions.

\begin{table}[tbp]
\centering\small
\caption{The buffer-reserve policy and \textsc{PARK}. Both charge for scarce
resources and account for the value of future cache reuse.}
\label{tab:br-park}
\begin{tabular}{@{}p{0.16\linewidth}p{0.37\linewidth}p{0.39\linewidth}@{}}
\toprule
Decision & Buffer-reserve policy (Algorithm~\ref{def:br-policy})
         & \textsc{PARK} (Algorithm~\ref{alg:online-dual})\\
\midrule
Resource prices &
Re-solving the buffered fluid LP each epoch yields shadow prices for
computation and memory. &
Updated online via projected subgradient steps driven by excess desired demand for work and active KV memory.\\
Routing &
The fluid bid accounts for current resource use and the change in
continuation value over the remaining horizon. &
The bid accounts for predicted prefix reuse now and the priced cost of one
subsequent warm turn after a migration.\\
Admission &
Use the buffered plan to route arrivals, rejecting an arrival when its
routes have no remaining quota. &
Offer an eligible endpoint with sufficient residual capacity; requests
without a qualifying offer wait for reconsideration.\\
Retention &
At each epoch boundary, weigh the continuation value of an idle cache
against the memory charge for keeping it. &
When active work needs memory, reclaim the lowest-priority blocks, where
priorities reflect reuse value per byte-epoch.\\
\bottomrule
\end{tabular}
\end{table}

\paragraph{Admission and routing.}
For a ready request $i$ and endpoint $e$, let $w_{ie}$ denote the estimated computation work for the uncached prompt and generated reply, and let $m_{ie}$ denote the estimated active KV-memory footprint. For a returning session, the controller credits expected prefix reuse only at its current cache owner $o_i$. It assigns endpoint $e$ the bid score
\begin{equation}
 s_{ie}=1+\omega_i-\lambda_e w_{ie}-\gamma_e m_{ie}+\Delta_{ie}.
 \label{eq:exp-online-dual-score}
\end{equation}
Here, $1$ represents the unit reward for serving the turn, $\omega_i$ is a bounded bonus that increases with waiting time prior to admission, and $\lambda_e$ and $\gamma_e$ price computation work and active memory, respectively. For a new session or a return to the current owner, $\Delta_{ie}=0$. For a migration away from the owner, $\Delta_{ie}=\alpha(F_{io_i}-F_{ie})$, where $F_{ie}$ is the priced resource cost of serving one subsequent warm turn at endpoint $e$ and $\alpha$ is a lookahead weight. Thus, migrating a session incurs the immediate cost of rebuilding its prefix, but earns credit for lower congestion on the next turn.

The controller restricts eligibility to endpoints with recent status reports. A returning session may migrate only when the candidate destination exhibits a sustained load advantage, offers a sufficiently higher bid than the current owner, and has elapsed a minimum cooldown since its previous endpoint binding. Standard requests are processed in decreasing order of their best eligible bid: each receives an offer for its highest-bidding eligible route that has a nonnegative bid and sufficient residual capacity after accounting for ongoing work and outstanding offers, while unplaced requests remain held in queue. Overdue requests are processed first in arrival order and may bypass the bid nonnegativity check, subject to the same capacity limits.

\paragraph{Updating prices.}
For each resource $j$ (computation work or active memory), let $C_{ej}$ denote the capacity budget at endpoint $e$ and $D_{ej}$ its desired load, comprising ongoing work, outstanding offers, and the demand of every ready request whose best eligible bid is nonnegative and selects $e$, which is measured \emph{before} capacity checks. Writing $p_{ej}$ for the unit shadow price ($\lambda_e$ or $\gamma_e$) and $z_{ej}=C_{ej}p_{ej}$ for its budget-scaled value, the controller updates prices with step size $\eta>0$ via
\begin{equation}
 z_{ej}^{+}=\bigl[z_{ej}+\eta(D_{ej}/C_{ej}-1)\bigr]_+ .
 \label{eq:exp-online-dual-update}
\end{equation}
Measuring desired demand prior to capacity filtering allows an endpoint's shadow price to keep rising even after its admissions are capped.

\paragraph{KV cache retention and eviction.}
For an idle prefix of length $L$ and idle age $a$, let $b_s=\beta L$ denote its memory footprint and $\widehat{\deff}(a,L)$ the estimated holding time per successful reuse, as defined in~\eqref{eq:bp-effective-return}. Each endpoint evaluates the retention index
\begin{equation}
 \widehat I(L,a)=\frac{\widehat g(L,a)}{b_s\,\widehat{\deff}(a,L)},
 \qquad
 \widehat g(L,a)=\min\{\widehat R,\,\lambda_e(a_1L+a_2L^2)\}.
 \label{eq:bp-retention-index}
\end{equation}
The numerator prices the work saved upon reuse at the estimated execution shadow price $\lambda_e$, using linear and quadratic rebuild-cost coefficients $a_1$ and $a_2$, and caps the resulting benefit at a fixed reference value; the index is set to zero if no useful return is forecast. Because the memory allocator reclaims physical KV blocks rather than whole sessions, it quantizes each session's index into a discrete priority and sums those priorities across all sessions sharing a given block. Unpinned, reclaimable blocks are bucketed by aggregate priority and evicted from the lowest nonempty bucket first, allowing shared prefix blocks to survive eviction even when individual sharing sessions have modest standalone priority.

\begin{algorithm}[tbp]
\caption{\textsc{PARK} (Price-based Admission and Routing with KV-cache Retention)}
\label{alg:online-dual}\label{alg:park}
\small
\emph{Central controller, each control round.}
\begin{enumerate}[leftmargin=*,itemsep=2pt,topsep=3pt]
\item Refresh endpoint loads and outstanding offers, and estimate the work
and active memory for each ready request at each eligible destination.
\item Compute route bids using~\eqref{eq:exp-online-dual-score}, and record
desired load at each request's highest-bidding eligible destination whenever
that bid is nonnegative, prior to capacity filtering.
\item Process overdue requests in arrival order, followed by standard requests
in decreasing order of best eligible bid; issue an offer for each request's
highest-bidding qualifying route, reserve its capacity, and hold unplaced
requests.
\item Update work and memory prices using~\eqref{eq:exp-online-dual-update}.
When a client acts on an offer, recheck eligibility and capacity before
dispatch, and release unused offers upon expiry.
\end{enumerate}
\emph{Endpoint allocator, whenever active work needs KV blocks.}
Reclaim unpinned blocks from the lowest nonempty priority bucket until
sufficient blocks are available or none remain eligible.
\end{algorithm}

}
\section{Experimental Results}
\label{sec:experiments}

We evaluate our framework using real multi-turn AI-agent workloads on large-scale GPU clusters. Section~\ref{sec:exp-setup} describes the experimental environment, workload characteristics, benchmark policies, and performance metrics. Section~\ref{sec:exp-main} compares \textsc{PARK} (Algorithm~\ref{alg:online-dual}) against standard routing and caching baselines across 16- and 128-GPU deployments. Finally, Section~\ref{sec:exp-mechanisms} isolates the two mechanisms driving these gains: value-aware cache retention and session-affine routing. 

\subsection{Experimental Setting, Workloads, and Benchmarks}
\label{sec:exp-setup}

\paragraph{Agent execution and cluster architecture.}
In each experiment, a bounded-concurrency client manages a pool of \emph{sessions}, each corresponding to one execution, or \emph{rollout}, of an agent task. Within each session, the agent alternates between LLM inference calls and external tool execution, such as running shell commands or tests in containerized environments. Each admitted inference \emph{request} initiates a turn, and all requests within a session carry the same persistent identifier.

All cluster experiments serve the same fixed Step~3.7~Flash checkpoint~\citep{stepfun2026step37flash}, a mixture-of-experts (MoE) model with approximately 197~billion parameters, 45~decoder layers, and hidden dimension 4,096 (its 42~MoE layers each activate eight of 288~routed experts per token alongside one shared expert). Each serving replica spans eight H200 GPUs using tensor parallelism (TP8). Consequently, the 8-, 16-, and 128-GPU deployments provide one, two, and sixteen independently routable endpoints, respectively. No model weights are updated during measurement. A central CPU router queries the request-scheduler over gRPC for admission and endpoint selection, while vLLM manages token-level batching and physical KV blocks locally on each replica. Additional telemetry, instrumentation, and hyperparameter details appear in Appendix~\ref{app:exp-details}.

\paragraph{Workloads.}
We evaluate two multi-turn software-engineering benchmarks:
\begin{itemize}
  \item \textbf{SWE:} A production mixture of 2,645 software-engineering tasks drawn from SWE-Factory (1,632)~\citep{guo2026swefactory}, SWE-rebench-v2 (975)~\citep{badertdinov2026swebenchv2}, and Multi-SWE-RL (38)~\citep{zan2025multiswebench,multiswerl2025dataset}, spanning 2,641 distinct task identifiers across multiple programming languages.
  \item \textbf{TB:} A corpus of 1,468 distinct hard software-engineering tasks packaged in the Terminal-Bench~2.0 task format~\citep{merrill2026terminalbench} and run on the Terminus scaffold~\citep[Section~3]{merrill2026terminalbench}.
\end{itemize}

Table~\ref{tab:exp-trajectories} and Figure~\ref{fig:exp-trajectories} summarize the two workloads in the three-hour \textsc{PARK} runs used for Table~\ref{tab:exp-routing}. Although both workloads accumulate long prompt histories (averaging 43,168 prompt tokens per request on SWE and 36,254 on TB), their interaction patterns differ sharply. Completed SWE sessions average 71.87 turns with relatively short responses (mean 433.46 output tokens), creating frequent opportunities to reuse the conversation's prefix across rapid tool calls. By contrast, completed TB sessions average 16.52 turns with responses approximately 5.3 times as long on average (2,290.16 output tokens), placing substantially heavier pressure on active decode concurrency and KV-cache capacity per turn.

\noindent
\begin{table}
\begin{minipage}[b]{0.45\textwidth}
\centering\small
\resizebox{\textwidth}{!}{
\begin{tabular}{lrr}
\toprule
Statistic & SWE & TB \\
\midrule
Total sessions & 2,244 & 2,521 \\
Completed sessions & 1,769 & 2,045 \\
Completed turns & 162,894 & 39,667 \\
Prompt tokens per turn & 40,841 / 95,359 & 35,373 / 91,361 \\
Output tokens per turn & 147 / 1,014 & 808 / 5,589 \\
\bottomrule
\end{tabular}}
\medskip
\caption{Workload trajectories in the 16-GPU \textsc{PARK} runs of Table~\ref{tab:exp-routing}. Paired entries give the median / 90th percentile.} \label{tab:exp-trajectories}
\end{minipage}\hfill \begin{minipage}[b]{0.55\textwidth}
\centering\small
\includegraphics[width=\linewidth]{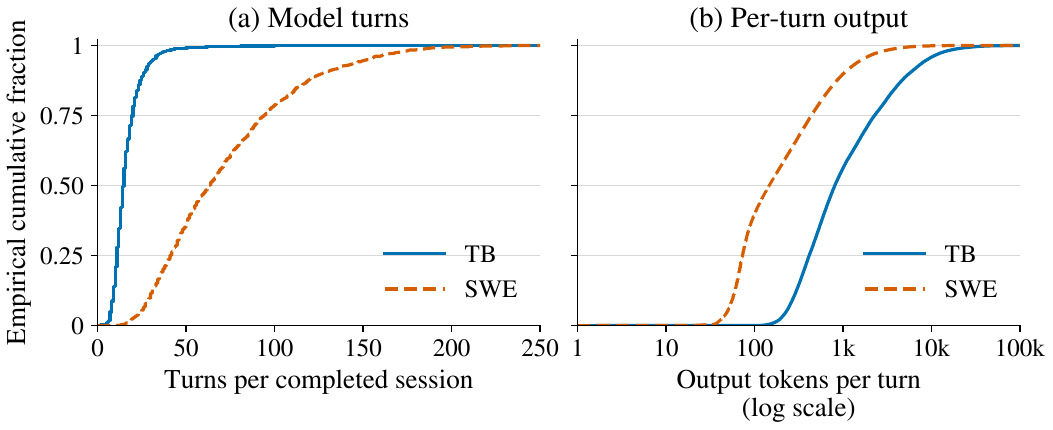}
\captionof{figure}{Empirical cumulative distributions of (left) model turns per completed session, and (right) output tokens per turn on SWE and TB.}
\label{fig:exp-trajectories}
\end{minipage}
\end{table}

\paragraph{Benchmark policies.}
Table~\ref{tab:exp-policies} summarizes the six serving policies compared in our cluster experiments. All five baseline policies admit every arriving prompt immediately without holding and rely on vLLM's native least-recently-used (LRU) eviction when an endpoint needs KV blocks; they differ only in how they select an endpoint. Let $R_e$ and $W_e$ denote running and waiting requests at endpoint $e$, $B_e$ its reported pending-prefill tokens, and $U_e = P - \min\{H_e, P-1\}$ the estimated uncached prefill tokens for a prompt of length $P$ given longest witnessed local prefix $H_e$. The three queue- and prefill-based baselines minimize $R_e+W_e$ (\emph{Shortest queue}), $R_e+4W_e$ (\emph{Weighted queue}), and $B_e$ (\emph{Least prefill}), respectively. \emph{SGLang routing} ports SGLang's \texttt{cache\_aware\_zmq} endpoint-selection  rule~\cite{sglang-cache-aware} into our shared vLLM serving stack: it routes to the endpoint with the longest cached prefix match whenever the relative match exceeds $50\%$ and server load imbalance remains within configured thresholds, and otherwise selects the least-loaded endpoint. In contrast, \textsc{PARK} (Algorithm~\ref{alg:online-dual}) couples shadow-price admission and routing with priority-bucket KV eviction. Exact tie-breaking and staleness rules appear in Appendix~\ref{app:exp-details}.

\begin{table}[ht]
\centering\small
\caption{Summary of evaluated policies.}
\label{tab:exp-policies}
\begin{tabular}{lll}
\toprule
Policy & Admission and Routing Rule & KV Eviction Rule \\
\midrule
Round robin & cyclic rotation across replicas & Native LRU \\
Shortest queue & $\arg\min_e (R_e + W_e)$ & Native LRU \\
Weighted queue &  $\arg\min_e (R_e + 4W_e)$ & Native LRU \\
Least prefill & $\arg\min_e B_e$ & Native LRU \\
SGLang routing (port) & prefix-match threshold with load fallback & Native LRU \\
\textsc{PARK} (Algorithm~\ref{alg:online-dual}) & Bid-price session-affine routing & Priority-bucket index \\
\bottomrule
\end{tabular}
\end{table}

\paragraph{Performance metrics.}
Our primary serving metric is \emph{output throughput} (generated tokens/s), defined as the aggregate rate at which output tokens are produced across all replicas—computed from each replica's cumulative generated-token counter divided by elapsed wall-clock time, as recorded by an independent monitor (Appendix~\ref{app:exp-details}). Because output throughput tracks continuous token generation while requests are decoding, whereas \emph{turn throughput} (completed turns/min) registers only when a turn finishes, the two rates can diverge over a run: synchronized client release rounds or short tool-call replies produce sharp early spikes in completed turns/min with modest token/s, whereas long generations sustain high token/s at lower turn-completion rates. We therefore report output tokens/s alongside completed turns/min and completed sessions.

\subsection{\textsc{PARK} Performance on Real Agent Workloads}
\label{sec:exp-main}

\begin{figure}[tbp]
 \centering
 \includegraphics[width=\linewidth]{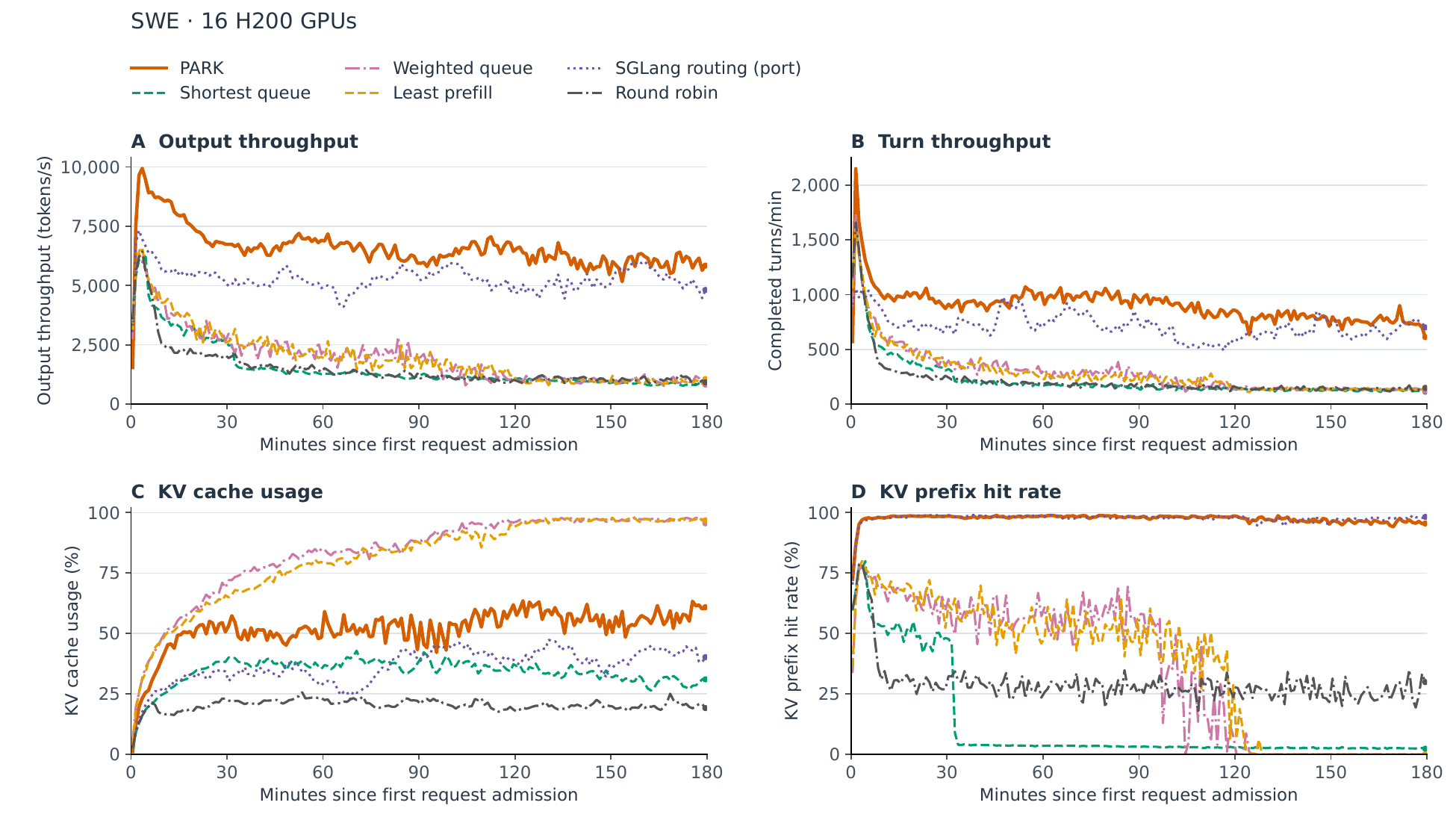}
 \caption{Time-series trajectories of output throughput, turn throughput, KV cache occupancy, and KV prefix hit rate on SWE (16~GPUs).}
 \label{fig:swe16-policy-timeseries}
\end{figure}

\begin{figure}[tbp]
 \centering
 \includegraphics[width=\linewidth]{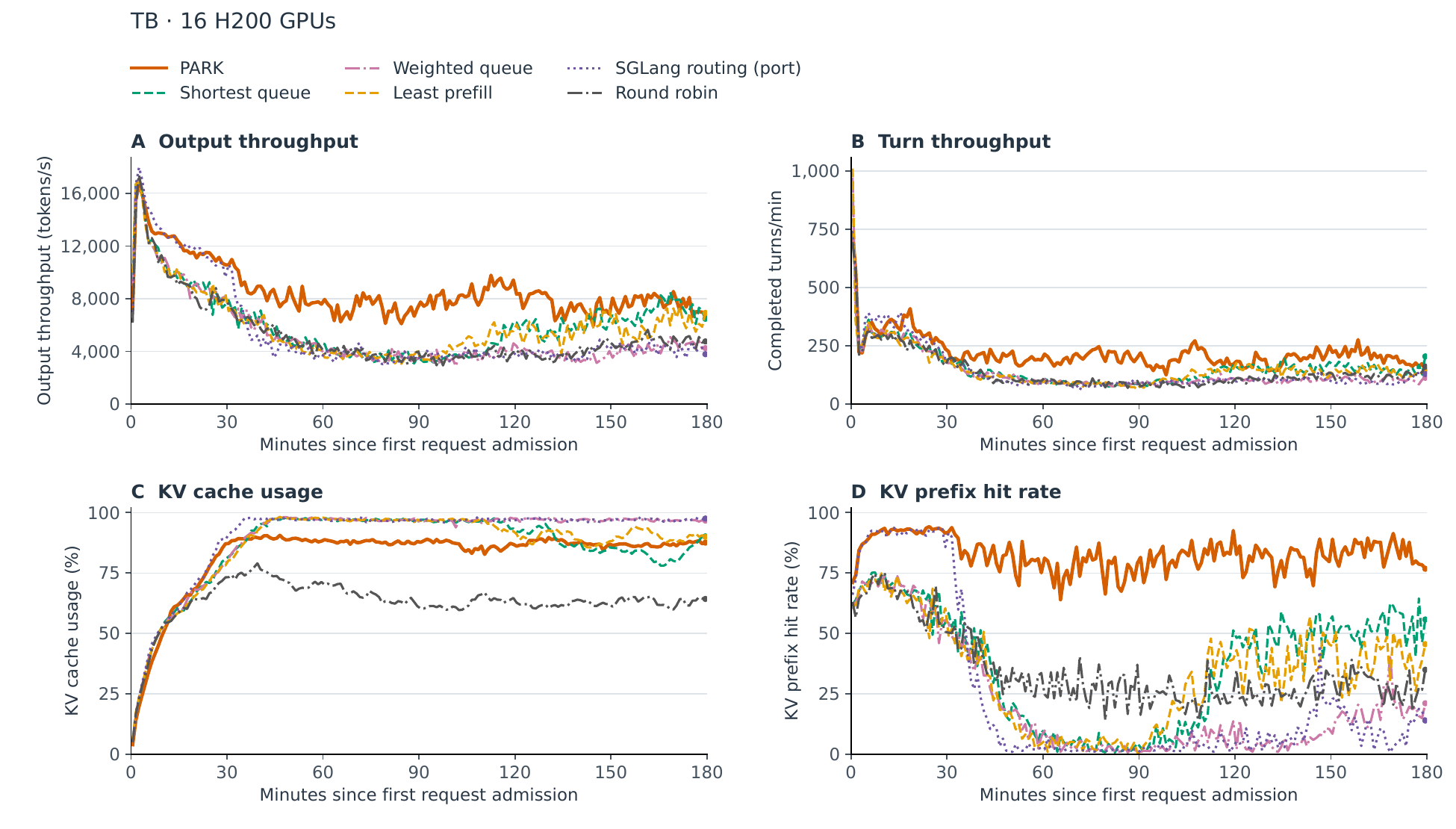}
 \caption{Time-series trajectories of output throughput, turn throughput, KV cache occupancy, and KV prefix hit rate on TB (16~GPUs).}
 \label{fig:tb16-policy-timeseries}
\end{figure}

\paragraph{16-GPU policy comparison.}
Table~\ref{tab:exp-routing} reports full-window performance across all six policies on the two-replica (16-GPU) deployment (three hours of measurement), and Figures~\ref{fig:swe16-policy-timeseries} and~\ref{fig:tb16-policy-timeseries} show the corresponding time-series trajectories on SWE and TB. On both workloads, \textsc{PARK} achieves the highest output throughput, delivering 6,575~tokens/s on SWE and 8,590~tokens/s on TB.

\begin{table}[ht]
\centering\small
\caption{Full-window performance on 16~H200 GPUs (two TP8 replicas). Each three-hour window starts at the first request admission.}
\label{tab:exp-routing}
\begin{tabular}{lrrrrrr}
\toprule
 & \multicolumn{3}{c}{SWE} & \multicolumn{3}{c}{TB} \\
\cmidrule(lr){2-4} \cmidrule(lr){5-7}
Policy & Tokens/s & vs.\ Best Base & Turns/min & Tokens/s & vs.\ Best Base & Turns/min \\
\midrule
Round robin & 1,505 & $-71.6\%$ & 215.8 & 5,261 & $-14.2\%$ & 136.9 \\
Shortest queue & 1,551 & $-70.7\%$ & 221.7 & 6,135 & $0.0\%$ & 156.9 \\
Weighted queue & 1,946 & $-63.3\%$ & 290.3 & 5,204 & $-15.2\%$ & 136.7 \\
Least prefill & 1,959 & $-63.0\%$ & 283.3 & 5,920 & $-3.5\%$ & 151.1 \\
SGLang routing (port) & 5,296 & $0.0\%$ & 716.9 & 5,509 & $-10.2\%$ & 140.6 \\
\textsc{PARK} & \textbf{6,575} & $\mathbf{+24.2\%}$ & \textbf{905.0} & \textbf{8,590} & $\mathbf{+40.0\%}$ & \textbf{220.4} \\
\bottomrule
\end{tabular}
\end{table}

Over the full three-hour window (Table~\ref{tab:exp-routing}), \textsc{PARK} outperforms the strongest measured baseline on each workload. On SWE, it reaches 6,575~tokens/s, $24.2\%$ above \emph{SGLang routing} and $336.9\%$ above \emph{Round robin}. On TB, it reaches 8,590~tokens/s, $40.0\%$ above \emph{Shortest queue} and $63.3\%$ above \emph{Round robin}. The identity of the strongest baseline changes with the workload, consistent with the trajectory differences in Table~\ref{tab:exp-trajectories}. SWE sessions average about 72 turns, with short responses and repeated opportunities to reuse a conversation prefix. On SWE, \emph{SGLang routing} achieves a $97.8\%$ prefix hit rate and 5,296~tokens/s, compared with 1,551--1,946~tokens/s for the two queue-based policies. On TB, where responses average 2,290 output tokens, \emph{Shortest queue} leads the baselines at 6,135~tokens/s, ahead of \emph{SGLang routing} at 5,509~tokens/s. Across these two regimes, \textsc{PARK}'s joint admission and routing policy also achieves the highest turn throughput, completing 905.0~turns/min on SWE and 220.4~turns/min on TB.

\begin{figure}[tbp]
 \centering
 \includegraphics[width=\linewidth]{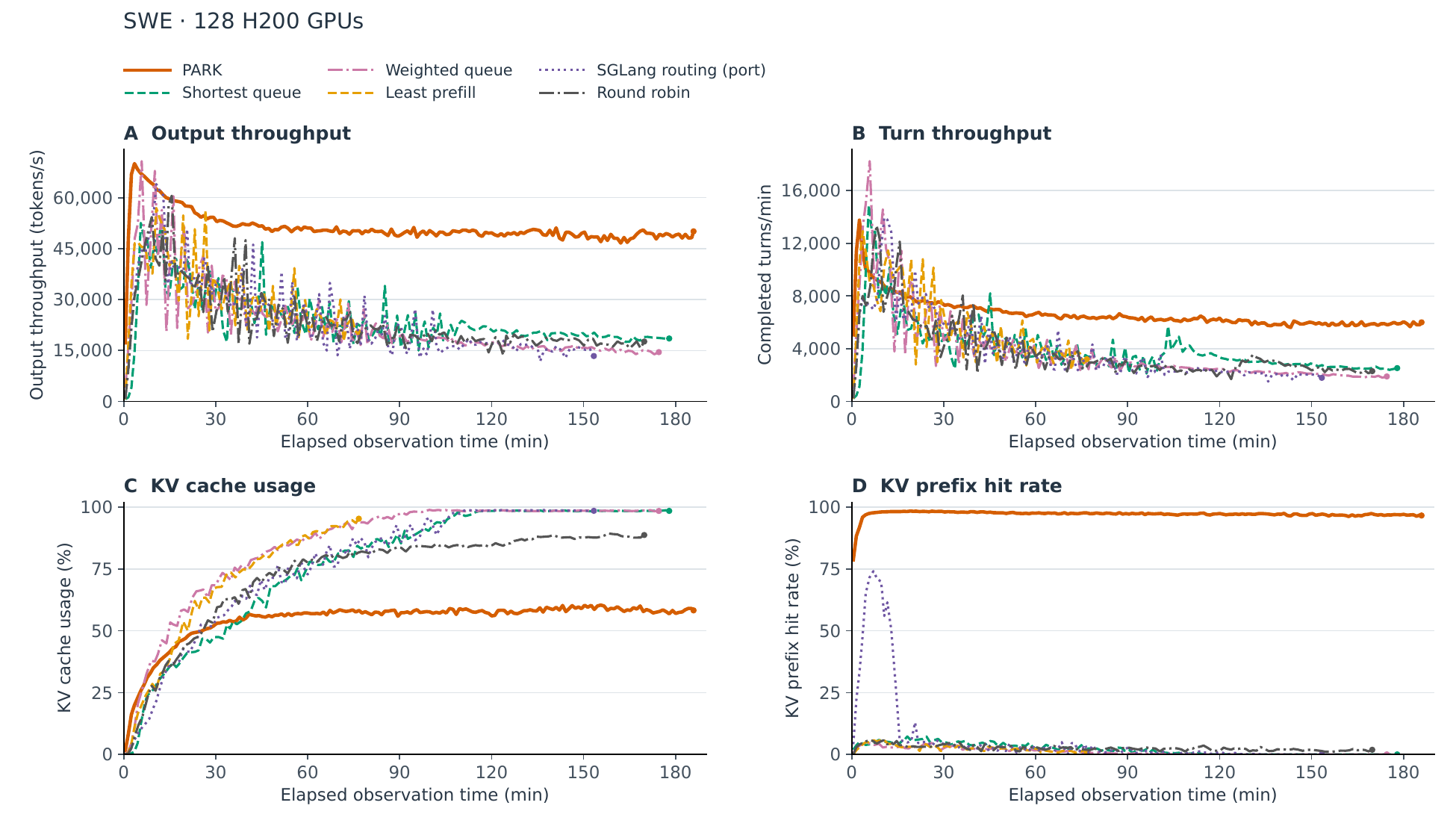}
 \caption{Time-series trajectories of output throughput, turn throughput, KV cache occupancy, and KV prefix hit rate on SWE (128~GPUs).}
 \label{fig:swe128-policy-timeseries}
\end{figure}

\paragraph{128-GPU deployment.}
To verify that \textsc{PARK}'s centralized primal--dual control and local eviction coordinate larger clusters under real agent traffic, we scale the deployment eightfold from two to sixteen TP8 replicas (128~GPUs) and
increase the maximum number of concurrent sessions eightfold from 512 to 4,096 (Table~\ref{tab:swe-scaleup}).
Over a three-hour measurement window, \textsc{PARK} sustains approximately 50,789~output tokens/s on SWE, substantially
outperforming the baselines (Figure~\ref{fig:swe128-policy-timeseries}).
The cluster completes 6,473~turns/min with a $97.42\%$ prefix cache hit rate and an average KV occupancy of $56.15\%$, and during this window, 12,295~sessions finished.
In a separate three-hour evaluation on TB at the same scale, \textsc{PARK} achieves 67,363~output tokens/s and 1,575~turns/min, with 13,563~completed sessions, while maintaining
a $91.32\%$ prefix cache hit rate at an average KV occupancy of $41.04\%$.

\begin{table}[tbp]
 \centering\small
 \caption{Cluster configurations for the \textsc{PARK} scale-up evaluation.}
 \label{tab:swe-scaleup}
 \begin{tabular}{lrr}
 \toprule
 & 16-GPU Deployment & 128-GPU Deployment \\
 \midrule
 Inference GPUs & 16 & 128 \\
 Routable TP8 replicas & 2 & 16 \\
 GPUs per replica & 8 & 8 \\
 Maximum pending sessions & 512 & 4,096 \\
 Permitted sessions per replica & 256 & 256 \\
 \bottomrule
 \end{tabular}
\end{table}

\subsection{Cache Retention and Routing Behavior}
\label{sec:exp-mechanisms}

\paragraph{Single-replica cache-retention ablation.}
To isolate the impact of index-based eviction from routing, we evaluate a single eight-GPU replica ($E=1$) on SWE with fixed round-robin dispatch and a maximum of 256 pending sessions (five-minute warmup followed by 30~minutes of measurement). Both arms use identical vLLM asynchronous scheduling and chunked prefill; only the selection of reclaimable idle KV blocks differs between native LRU and our priority-bucket retention index (Section~\ref{sec:swe-gradient-implementation}).

\begin{table}[ht]
\centering\small
\caption{Single-replica (8-GPU) cache-retention ablation on SWE (30-minute measurement window).}
\label{tab:exp-retention}
\begin{tabular}{lrr}
\toprule
Metric & Native LRU & Retention Index \\
\midrule
\multicolumn{3}{l}{\emph{Cache preservation}} \\
Prior-prefix token loss (\%) & 13.26 & \textbf{2.17} \\
Returning turns experiencing prefix loss (\%) & 9.48 & \textbf{2.72} \\
Verified prefix rebuild (M tokens) & 30.77 & \textbf{6.62} \\
\midrule
\multicolumn{3}{l}{\emph{Prefill computation}} \\
Total scheduled prefill (M tokens) & 55.52 & \textbf{29.70} \\
Local prefill tokens per output token & 13.26 & \textbf{5.19} \\
\midrule
\multicolumn{3}{l}{\emph{Throughput and latency}} \\
Output throughput (tokens/s) & 2,326 & \textbf{3,179}  \\
Completed model turns/min & 275.6 & \textbf{373.5} \\
Mean TTFT (s) & 25.56 & \textbf{16.49} \\
Mean ITL (ms) & 50.54 & \textbf{39.90} \\
Mean request latency (s) & 48.09 & \textbf{36.00} \\
\bottomrule
\end{tabular}
\end{table}

Table~\ref{tab:exp-retention} traces the causal chain from cache survival to serving throughput:
\begin{enumerate}
  \item \emph{Reduced prior-prefix loss:} By prioritizing intermediate-length prefixes with high expected reuse value per byte-epoch rather than evicting purely by recency, the retention index cuts prior-prefix token loss—the token-weighted share of previously verified prefix tokens missing when a session returns—from $13.26\%$ to $2.17\%$, reducing verified prefix recomputation from 30.77M to 6.62M tokens.
  \item \emph{Lower prefill overhead:} Preserving resident prefixes cuts total scheduled prefill nearly in half (from 55.52M to 29.70M tokens), reducing prefill tokens processed per generated output token from 13.26 to 5.19.
  \item \emph{Higher throughput and lower latency:} Freeing compute capacity from redundant prefill raises output throughput by $36.7\%$ (from 2,326 to 3,179~tokens/s) and reduces mean time-to-first-token (TTFT) by $35.5\%$ (from 25.56\,s to 16.49\,s).
\end{enumerate}

\paragraph{Multi-replica endpoint affinity and cache thrashing.}
In multi-replica clusters, routing and retention reinforce one another: returning a session to its cache owner creates an opportunity for prefix reuse, which generates savings only if the prefix survives the inter-turn tool delay. In the three-hour, 16-GPU SWE runs, \textsc{PARK} switched replicas on $0.14\%$ of returning requests (222 of 161,218), compared with $50.14\%$ for \emph{Round robin} (19,310 of 38,515).
Panels~C and~D of Figures~\ref{fig:swe16-policy-timeseries} and~\ref{fig:tb16-policy-timeseries} illustrate the operational consequence across time: baselines such as \emph{Shortest queue} accumulate high KV occupancy (Panel~C) while suffering near-zero prefix hit rates (Panel~D) because cross-replica switching and recency eviction continually thrash cached prefixes, whereas \textsc{PARK} maintains moderate KV occupancy ($\approx 50\%$ on SWE) alongside a $\approx 98\%$ prefix hit rate.

\section{Conclusion}
\label{sec:conclusion}

Multi-turn AI agents transform LLM serving from stateless request dispatch into a dynamic resource-allocation problem with endogenous state: routing determines where a conversation's prefix can next be reused, whereas cache retention determines whether that prefix survives inter-turn delays. By formulating a fluid approximation of this joint system, we characterized optimal admission, routing, and eviction through endpoint-specific shadow prices for computation and memory coupled with session continuation values (Theorem~\ref{prop:bp-control}). Building on this decomposition, we developed a synchronous buffer-reserve policy (Algorithm~\ref{def:br-policy}) and proved that setting aside a capacity reserve of order $\sqrt{n\log n}$ keeps emergency interventions rare and achieves an $O(\sqrt{n\log n})$ completion shortfall relative to the fluid benchmark (Theorem~\ref{thm:br-regret}).

To bring these structural principles to production clusters characterized by asynchronous arrivals, delayed load reports, and engine-local memory management, we designed \textsc{PARK} (Algorithm~\ref{alg:online-dual}). \textsc{PARK} pairs projected subgradient price updates and one-turn migration lookahead at a central controller with a priority-bucket KV retention index based on reuse value per byte-epoch at each endpoint. Across 16- and 128-GPU deployments serving multi-turn software-engineering workloads (Section~\ref{sec:experiments}), \textsc{PARK} outperforms both load-balancing and prefix-aware baselines across contrasting interaction regimes, adapting between long-horizon, short-reply trajectories (SWE) and memory-intensive, long-generation trajectories (TB). The single-replica ablation demonstrates the benefit of value-aware cache retention, while the multi-replica measurements document PARK’s high prefix-reuse rates.

\bibliographystyle{informs2014}
\bibliography{references}

\ECSwitch
\ECHead{Electronic Companion}

\renewcommand*{\theHsection}{EC.\arabic{section}}

\section{Proofs for bid-price control and retention}
\label{app:bp-proofs}

\subsection{Single-session values and optimal resource prices}
\label{app:bp-dual}

Fix a feasible remaining-horizon LP \eqref{eq:fluid-lp-t},
starting at epoch $t_0$ with population $m^0$. Write
$x=(m,\widetilde m,y,z)$ for its flows and $F(x)$ for its
objective. All state sums below are over the finite live state
space. The transition kernels are sub-stochastic, with missing
mass representing departure.

The flow equations imply
\[
 \|m(\tau+1)\|_1=\|\widetilde m(\tau+1)\|_1
 \le \|m(\tau)\|_1+\sum_{s\in\cSI}\nu_s(\tau).
\]
Thus the feasible set is
compact, the optimum is finite and attained, and LP strong
duality gives an attained dual optimum. 

For any fixed nonnegative dual prices $(\lambda,\gamma)$, define
the session values backward from
$v_s(H)=-\1\{s\in\cSA\}$ at the terminal. At each epoch $\tau$, set
$w_s(\tau)$ according to \eqref{eq:bp-boundary-option} for
cached idle states, and set $w_s(\tau)=v_s(\tau+1)$ for active
and cache-free idle states. Put $r_{s,0}(\tau)=0$ and use
\eqref{eq:bp-resource-rent} for real endpoints. The value of
action $e$ for an arriving prompt is
\begin{equation}
 A_{s,e}(\tau)=\1\{e\ne0\}-r_{s,e}(\tau)
       +\sum_{s'\in\cS}P_\tau^e(s,s')w_{s'}(\tau),
 \qquad s\in\cSI,\quad e\in\cEzero.
 \label{eq:bp-action-value}
\end{equation}
Its maximum $\eta_s(\tau)$ is conditional on a prompt arriving.
The beginning-of-epoch values are
\begin{align}
 \eta_s(\tau)&=\max_{e\in\cEzero}A_{s,e}(\tau),
       &&s\in\cSI,\nonumber\\
 v_s(\tau)&=p_s(\tau)\eta_s(\tau)
       +(1-p_s(\tau))\sum_{s'}Q_\tau^I(s,s')w_{s'}(\tau),
       &&s\in\cSI,\label{eq:bp-canonical-idle}\\
 v_s(\tau)&=-\sum_e\left\{\lambda_e(\tau)\ell_{s,e}(\tau)
                  +\gamma_e(\tau)\kappa_{s,e}^{\mathrm{act}}(\tau)\right\}
       +\sum_{s'}Q_\tau^A(s,s')w_{s'}(\tau),
       &&s\in\cSA.\label{eq:bp-canonical-active}
\end{align}
These equations define values at every state, including states
with zero flow in an optimal solution. The following lemma expresses the dual objective of the fluid linear program~\eqref{eq:fluid-lp-t} in terms of the execution and KV-memory shadow prices $(\lambda,\gamma)$, establishes strong duality, and decomposes the optimality gap $\mathcal D(\lambda,\gamma)-F(x)$ of any feasible flow into nonnegative capacity-slack and reduced-cost terms.

\begin{lemma}
\label{lem:bp-canonical}
The dual objective at $t_0$ with initial population $m^0$ has value
\begin{equation}
\begin{aligned}
 \mathcal D(\lambda,\gamma)
 ={}&\sum_{\tau=t_0}^{H-1}\sum_e
       \{\lambda_e(\tau)B_e+\gamma_e(\tau)G_e\}
 +\sum_s m_s^0v_s(t_0)
       +\sum_{\tau=t_0}^{H-1}\sum_{s\in\cSI}
                         \nu_s(\tau)\eta_s(\tau).
\end{aligned}
\label{eq:bp-dual-function}
\end{equation}
At optimal dual variable $(\lambda,\gamma)$, strong duality holds, i.e.,
\(
  \mathcal D(\lambda,\gamma)=\mathcal V_{t_0}(m^0).
\)
Moreover, for every feasible primal flow $x$ and nonnegative dual $(\lambda,\gamma)$, we have
\begin{align}
 \mathcal D(\lambda,\gamma)-F(x)
 ={}&
 \sum_{\tau,e}\lambda_e(\tau)
       \{B_e-\Phi_e^{\mathrm{exe}}(\tau;x)\}
 +\sum_{\tau,e}\gamma_e(\tau)
       \{G_e-\Phi_e^{\mathrm{KV}}(\tau;x)\}
 \nonumber\\
 &+\sum_{\tau,s\in\cSI,a\in\cEzero}
       y_{s,a}(\tau)\{\eta_s(\tau)-A_{s,a}(\tau)\}
 \nonumber\\
 &+\sum_{\tau,s\in\cSI:o(s)\ne0}z_{s,\KEEP}(\tau)
       \{w_s(\tau)-v_s(\tau+1)+\gamma_{o(s)}(\tau)b_s\}
 \nonumber\\
 &+\sum_{\tau,s\in\cSI}z_{s,\DROP}(\tau)
       \{w_s(\tau)-v_{\Ev(s)}(\tau+1)\}.
 \label{eq:bp-dual-gap}
\end{align}
Here $\Phi^{\mathrm{exe}},\Phi^{\mathrm{KV}}$ are the resource
left-hand sides in \eqref{eq:fluid-lp-t}.
Every summand is
nonnegative so $\mathcal D(\lambda,\gamma)\geq F(x)$.
\end{lemma}

\begin{proof}
Fix nonnegative dual prices $(\lambda,\gamma)$ and the associated
backward values. Let $\mathcal X$ contain the nonnegative
flows satisfying the initial condition and all flow equations,
with resource constraints removed. We first derive the gap
identity (\ref{eq:bp-dual-gap}) using only these equations, and then use it to
identify the dual function.

\paragraph{1. Account for one epoch.}
Fix $\tau$, suppress current-epoch arguments, and write
$m_s^+=m_s(\tau+1)$,
$\widetilde m_s=\widetilde m_s(\tau+1)$, and
$v_s^+=v_s(\tau+1)$.
Define the nonnegative differences between the best action
value and each available action's value:
\[
\begin{aligned}
 \epsilon_{s,a}&=\eta_s-A_{s,a},
       &&s\in\cSI,\quad a\in\cEzero,\\
 \epsilon_{s,\KEEP}&=w_s-v_s^++\gamma_{o(s)}b_s,
       &&s\in\cSI,\quad o(s)\ne0,\\
 \epsilon_{s,\DROP}&=w_s-v_{\Ev(s)}^+,
       &&s\in\cSI.
\end{aligned}
\]
Let $R_\tau$ be admissions, $\mathcal C_\tau$ the total
resource charge, and $\mathcal I_\tau$ its idle-memory part:
\[
\begin{aligned}
 R_\tau&=\sum_{s\in\cSI,e\in\cE}y_{s,e},\\
 \mathcal C_\tau
 &=\sum_e\{\lambda_e\Phi_e^{\mathrm{exe}}
                   +\gamma_e\Phi_e^{\mathrm{KV}}\},\\
 \mathcal I_\tau
 &=\sum_{\substack{s\in\cSI:\\o(s)\ne0}}
                 \gamma_{o(s)}b_s z_{s,\KEEP}.
\end{aligned}
\]

\emph{Routing and service.}
Arrival conservation (\ref{eq:fluid-routing-conservation}) and the definition of $A_{s,a}$  (\ref{eq:bp-action-value}) give
\begin{align*}
 \sum_{s\in\cSI}(p_sm_s+\nu_s)\eta_s
 &=\sum_{s\in\cSI,a\in\cEzero}y_{s,a}
                                  (A_{s,a}+\epsilon_{s,a})\\
 &=R_\tau-\sum_{s\in\cSI,e\in\cE}r_{s,e}y_{s,e}
   +\sum_r\left[\sum_{s\in\cSI,a\in\cEzero}
                    y_{s,a}P_\tau^a(s,r)\right]w_r
   +\sum_{s\in\cSI,a\in\cEzero}y_{s,a}\epsilon_{s,a}.
\end{align*}
Meanwhile, weighting the idle and active value recursions
by their populations (i.e., multiplying (\ref{eq:bp-canonical-idle}) and (\ref{eq:bp-canonical-active}) with $m_s$ and $\nu_s$) yields
\begin{align*}
 \sum_s m_s v_s+\sum_{s\in\cSI}\nu_s\eta_s
 ={}&\sum_{s\in\cSI}(p_sm_s+\nu_s)\eta_s
 -\sum_{s\in\cSA}m_s\sum_e
       \{\lambda_e\ell_{s,e}
                   +\gamma_e\kappa_{s,e}^{\mathrm{act}}\}\\
 &+\sum_r\left[
       \sum_{s\in\cSA}m_sQ_\tau^A(s,r)
       +\sum_{s\in\cSI}(1-p_s)m_sQ_\tau^I(s,r)
          \right]w_r.
\end{align*}
Substitute the first display into the second. The coefficient
of $w_r$ becomes the pre-retention inflow
$\widetilde m_r$. Since
$r_{s,e}=\lambda_ec_{s,e}
          +\gamma_e\kappa_{s,e}^{\mathrm{adm}}$,
the resource terms sum to $-\mathcal C_\tau+\mathcal I_\tau$:
they include execution and unfinished-turn memory, with
idle-memory charges still excluded. Hence
\begin{equation}
 \sum_s m_s v_s+\sum_{s\in\cSI}\nu_s\eta_s
 =R_\tau-\mathcal C_\tau+\mathcal I_\tau
       +\sum_r\widetilde m_r w_r
       +\sum_{s\in\cSI,a\in\cEzero}y_{s,a}\epsilon_{s,a}.
 \label{eq:bp-telescope-first}
\end{equation}

\emph{Boundary retention.}
For a cached idle state, retention conservation (\ref{eq:fluid-retention-conservation}) gives
\[
 \widetilde m_s w_s
 =z_{s,\KEEP}
       (v_s^+-\gamma_{o(s)}b_s+\epsilon_{s,\KEEP})
  +z_{s,\DROP}(v_{\Ev(s)}^++\epsilon_{s,\DROP}).
\]
For a cache-free state, KEEP is zero and this equality holds
with only the DROP term, whose slack is zero.
The post-retention equation groups each DROP flow at its
destination $\Ev(s)$, so
\begin{align*}
 \sum_{r\in\cSI}m_r^+v_r^+
 &=\sum_{r\in\cSI}
      \left[z_{r,\KEEP}
          +\sum_{\substack{s\in\cSI:\\\Ev(s)=r}}
                       z_{s,\DROP}\right]v_r^+
 =\sum_{s\in\cSI}
        \{z_{s,\KEEP}v_s^+
                    +z_{s,\DROP}v_{\Ev(s)}^+\}.
\end{align*}
Active states satisfy $\widetilde m_s=m_s^+$ and $w_s=v_s^+$.
Summing their contributions and the idle-state identities gives
\begin{align}
 \sum_r\widetilde m_r w_r
 ={}&\sum_r m_r^+v_r^+-\mathcal I_\tau
       +\sum_{\substack{s\in\cSI:\\o(s)\ne0}}
                          z_{s,\KEEP}\epsilon_{s,\KEEP}
       +\sum_{s\in\cSI}z_{s,\DROP}\epsilon_{s,\DROP}.
 \label{eq:bp-telescope-second}
\end{align}

\paragraph{2. Sum over the horizon.}
Let $\mathcal E_\tau$ collect the routing and retention
value losses in this epoch:
\[
 \mathcal E_\tau
 =\sum_{s\in\cSI,a\in\cEzero}y_{s,a}\epsilon_{s,a}
  +\sum_{\substack{s\in\cSI:\\o(s)\ne0}}
                        z_{s,\KEEP}\epsilon_{s,\KEEP}
  +\sum_{s\in\cSI}z_{s,\DROP}\epsilon_{s,\DROP}.
\]
Combining \eqref{eq:bp-telescope-first} and
\eqref{eq:bp-telescope-second} cancels $\mathcal I_\tau$:
\[
 \sum_s m_s(\tau)v_s(\tau)
 +\sum_{s\in\cSI}\nu_s(\tau)\eta_s(\tau)
 -\sum_s m_s(\tau+1)v_s(\tau+1)
 =R_\tau-\mathcal C_\tau+\mathcal E_\tau.
\]
Summing over $\tau=t_0,\ldots,H-1$ cancels the intermediate
population values. Using $m(t_0)=m^0$,
$v_s(H)=-\1\{s\in\cSA\}$, and
$F(x)=\sum_\tau R_\tau-\sum_{s\in\cSA}m_s(H)$ gives
\[
 \sum_s m_s^0v_s(t_0)
 +\sum_{\tau,s\in\cSI}\nu_s(\tau)\eta_s(\tau)
 =F(x)-\sum_\tau\mathcal C_\tau+\sum_\tau\mathcal E_\tau.
\]
Adding the capacity-price constants and rearranging yields
\[
 \mathcal D(\lambda,\gamma)-F(x)
 =\sum_{\tau,e}\left\{
       \lambda_e(\tau)[B_e-\Phi_e^{\mathrm{exe}}(\tau;x)]
      +\gamma_e(\tau)[G_e-\Phi_e^{\mathrm{KV}}(\tau;x)]
                      \right\}
       +\sum_\tau\mathcal E_\tau.
\]
Expanding $\mathcal E_\tau$ is exactly
\eqref{eq:bp-dual-gap}. For a capacity-feasible flow, every
term is nonnegative: resource slacks, prices, flows, and
action-value differences are all nonnegative.

\paragraph{3. Identify the dual function and its optimum.}
Define the resource Lagrangian by
\[
 \mathcal L(x;\lambda,\gamma)
 =F(x)+\sum_{\tau,e}\left\{
       \lambda_e(\tau)[B_e-\Phi_e^{\mathrm{exe}}(\tau;x)]
      +\gamma_e(\tau)[G_e-\Phi_e^{\mathrm{KV}}(\tau;x)]
                        \right\}.
\]
The calculation above used only flow equations. It therefore
applies to every $x\in\mathcal X$, and gives
\begin{equation}
 \mathcal D(\lambda,\gamma)-\mathcal L(x;\lambda,\gamma)
       =\sum_\tau\mathcal E_\tau\ge0.
 \label{eq:bp-relaxed-gap}
\end{equation}
This bound is attained. Starting from $m^0$, at every epoch
send each state's available arrival mass to an action
maximizing $A_{s,a}$. Propagate the resulting flows through
the transition equations, and assign each pre-retention idle
population to a maximizing KEEP or DROP action. Continue
through the finite horizon. The resulting nonnegative flows
satisfy every equation defining $\mathcal X$. Every positive
routing or retention flow has zero slack, so
$\mathcal E_\tau=0$ at every epoch. Consequently,
\(
 \max_{x\in\mathcal X}\mathcal L(x;\lambda,\gamma)
       =\mathcal D(\lambda,\gamma).
\)

Finally, take an optimal full LP dual solution with resource
prices $(\lambda^\star,\gamma^\star)$ and flow-equation
multipliers $\widehat d$. Write the flow equations as
$Ax=b$. By full dual feasibility and strong duality,
\[
 \mathcal L(x;\lambda^\star,\gamma^\star)
       +\widehat d^{\mathsf T}(b-Ax)
       \le\mathcal V_{t_0}(m^0),\qquad x\ge0.
\]
Indeed, the left side has nonpositive variable coefficients
and constant term equal to that dual's objective.
For $x\in\mathcal X$, $b-Ax=0$; maximizing therefore gives
$\mathcal D(\lambda^\star,\gamma^\star)
\le\mathcal V_{t_0}(m^0)$.
Conversely, any optimal capacity-feasible flow $x^\star$
has nonnegative resource slacks, so
\[
 \mathcal V_{t_0}(m^0)=F(x^\star)
 \le\mathcal L(x^\star;\lambda^\star,\gamma^\star)
 \le\mathcal D(\lambda^\star,\gamma^\star).
\]
Thus equality holds. In particular, pairing optimal prices
with any optimal primal flow makes every term in
\eqref{eq:bp-dual-gap} zero, as used in the bid-price theorem.
\end{proof}

\subsection{Proof of the bid-price characterization}
\label{app:bp-control}

\begin{proof}[Proof of Theorem~\ref{prop:bp-control}]
According to (\ref{eq:bp-dual-gap}) in
Lemma~\ref{lem:bp-canonical}, complementary slackness for
a routing variable gives
\(
 y_{s,e}(t)\{\eta_s(t)-A_{s,e}(t)\}=0.
\)
A positive flow therefore uses an action maximizing
$A_{s,e}(t)$. Since
\[
 A_{s,e}(t)-A_{s,0}(t)=\mathcal B_{s,e}(t)
 \qquad(e\in\cEzero),
\]
subtracting the common rejection value preserves the action
ordering and proves part (a).

For an owned idle state, complementary slackness gives
\begin{align*}
 z_{s,\KEEP}(t)
  \{w_s(t)-v_s(t+1)+\gamma_{o(s)}(t)b_s\}&=0,\\
 z_{s,\DROP}(t)
  \{w_s(t)-v_{\Ev(s)}(t+1)\}&=0.
\end{align*}
Write $g=\Delta_s(t+1)-\gamma_{o(s)}(t)b_s$, the KEEP
value minus the DROP value. If $g>0$, the boundary maximum
is KEEP, and the DROP slack equals $g$. Hence
$z_{s,\DROP}(t)=0$, and retention conservation implies
$z_{s,\KEEP}(t)=\widetilde m_s(t+1)$.
If $g<0$, the maximum is DROP, and the KEEP slack is $-g$.
Thus $z_{s,\KEEP}(t)=0$ and
$z_{s,\DROP}(t)=\widetilde m_s(t+1)$.
If $g=0$, both slacks vanish. Tied routing actions likewise
have zero slack, leaving their allocation to the flow and
capacity constraints. This proves part (b) and the tie statement.
\end{proof}

\subsection{Cache length, return time, and retention}
\label{app:bp-retention}

\subsubsection{The value available at the next prompt}
\label{app:bp-retention-deficit}

Suppose a finite return has positive probability. A policy
$\pi$ starts with the cache available at that return, observes
the date and prompt information, and is nonanticipative thereafter.
Let $N^\pi$ count its completions through $H$ and let
$\mathcal C^\pi$ contain its execution, active-memory, and
subsequent idle-memory charges. Pre-return waiting charges
are excluded. Define
\begin{equation}
 \varepsilon(L)=
 \inf_\pi\max_{\substack{1\le h\le N\\\Pp(T=h)>0}}
 \left\{R-\E^\pi[N^\pi-\mathcal C^\pi\mid T=h]\right\}.
 \label{eq:bp-retention-epsilon}
\end{equation}
This is the smallest worst-date shortfall from $R$; each
date's value averages over its prompt and subsequent service.

Let $s_h$ denote the observed cached state at a return in
epoch $t+h$. Starting from idle, admitted turns minus the
terminal active indicator equals completed turns on every
sample path. The admission rewards and terminal values in
the backward recursion therefore make $\eta_{s_h}(t+h)$
the optimal expected completions minus resource charges.
Backward induction selects a maximizing action at every
state and date, including every boundary-retention decision.
The resulting policy is optimal from every return state
simultaneously. Every other policy has weakly smaller value
at each state and hence at each return date. It follows that
\begin{equation}
 \varepsilon(L)=
 \max_{\substack{1\le h\le N\\\Pp(T=h)>0}}
 \left\{R-\E[\eta_{s_h}(t+h)\mid T=h]\right\}.
 \label{eq:bp-retention-epsilon-value}
\end{equation}
The bound $R$ and nonnegative charges give an upper value
bound of $R$. Rejecting every prompt and dropping every cache
gives value zero. Consequently $0\le\varepsilon(L)\le R$ and
\begin{equation}
 R-\varepsilon(L)
 \le\E[\eta_{s_h}(t+h)\mid T=h]\le R
 \label{eq:bp-retention-warm-value}
\end{equation}
for every possible return date.

\subsubsection{Proof of the retention and eviction conditions}
\label{app:bp-retention-proof}

\begin{proof}[Proof of Proposition~\ref{prop:bp-retention-length}]
Let us first briefly talk about the proof idea. We first bound the additional value of having
the cache when the next prompt arrives. We then express retention
as choosing how long to wait for that prompt, and subtract the
storage charge from the expected benefit. Retention is optimal if
some waiting deadline has positive gain; eviction is optimal if
every positive deadline has negative gain. The quantity $\deff$
allows us to obtain the stated conditions from these comparisons.

Fix the current cached state $s$ and write $s^0=\Ev(s)$.
All probabilities below condition on the current information,
and all values use the same fixed dual prices. If $N=0$
or no return is possible, keeping the cache pays $\gamma b_s>0$
without providing future value, so eviction is optimal.
Henceforth a finite return has positive probability. For brevity, put
\[
 C=\underline\lambda(a_1L+a_2L^2),
 \qquad D=\bar\lambda_e(a_1L+a_2L^2),
\]
where $C$ bounds the first cold turn's expected execution charge
from below, and $D$ bounds reuse savings from above under the
additional assumption for the short-cache comparison.

\paragraph{Step 1: Bound the benefit at the next prompt.}
For each $h$ with $\Pp(T=h)>0$, let $s_h$ be the cached state
observed when the prompt arrives in epoch $t+h$, and let
$s_h^0=\Ev(s_h)$ be its cache-free counterpart. Define
\begin{equation}
 g_h=\E[\eta_{s_h}(t+h)-\eta_{s_h^0}(t+h)\mid T=h].
 \label{eq:bp-retention-return-gain}
\end{equation}
Thus, $g_h$ is the extra optimal continuation value provided
by the cache at return. It includes later routing and retention
decisions, but excludes the storage charges paid before that return.

We first show that the cache-free value is at most $[R-C]_+$.
Consider any continuation policy $\mu$ starting without the cache
at this return. Let $N^\mu$ count its completed turns through $H$,
and let $\mathcal C^\mu$ be all its resource charges from the
return through $H$. Define the first admission epoch by
\[
 \sigma=\inf\{\tau\in\{t+h,\ldots,H-1\}:
                 \mu\text{ admits a prompt of this session at }\tau\},
 \qquad \inf\varnothing=\infty,
\]
and put $A_\tau=\{\sigma=\tau\}$. 

On $A_\tau$, suppose the first admitted turn belongs to session
$i$, is arrival $j$ in epoch $\tau$, and is assigned to endpoint
$f$. Using the execution requirements from Section~3, define
\begin{equation}
 \mathcal C_\tau^{\mathrm{first}}
 =\lambda_f(\tau)C_{\tau,j,f}
  +\sum_{\substack{u=\tau+1,\ldots,H-1:\
        \text{this turn is active at the start of }u}}
       \lambda_f(u)L_{u,i,f}.
 \label{eq:bp-retention-first-charge}
\end{equation}
This is the total execution charge of that first admitted turn
within the horizon. The sum stops when that turn completes;
it does not include later turns of the same session. Set
$\mathcal C_\tau^{\mathrm{first}}=0$ outside $A_\tau$.

Let $\mathcal H_\tau$ denote the information observed at the
routing decision, including the chosen action. On $A_\tau$,
there are no earlier completions, and the conditional bounds
on completions and cold execution charges give
\[
 \E^\mu[N^\mu\mid\mathcal H_\tau]\le R,
 \qquad
 \E^\mu[\mathcal C_\tau^{\mathrm{first}}
                         \mid\mathcal H_\tau]\ge C.
\]
The first-admission event $A_\tau$ is known at this decision,
as is the earlier return event $\{T=h\}$. Therefore, the
law of iterated expectations yields
\begin{align*}
 &\E^\mu[\1_{A_\tau}(N^\mu-\mathcal C_\tau^{\mathrm{first}})
                                                    \mid T=h]
=\E^\mu\!\left[
    \1_{A_\tau}\E^\mu[N^\mu-\mathcal C_\tau^{\mathrm{first}}
                           \mid\mathcal H_\tau]\,\middle|\,T=h\right]
\le\Pp^\mu(A_\tau\mid T=h)(R-C).
\end{align*}
On $A_\tau$, total charges satisfy
$\mathcal C^\mu\ge\mathcal C_\tau^{\mathrm{first}}$.
If $\sigma=\infty$, no turn completes and the value is zero.
The events $A_\tau$ are disjoint, so summing gives
\begin{align*}
 \E^\mu[N^\mu-\mathcal C^\mu\mid T=h]
 &\le\sum_{\tau=t+h}^{H-1}
       \Pp^\mu(A_\tau\mid T=h)(R-C)
 =\Pp^\mu(\sigma<\infty\mid T=h)(R-C)
 \le[R-C]_+.
\end{align*}
Optimizing over $\mu$, and observing that rejection of every
prompt gives value zero, proves
\begin{equation}
 0\le\E[\eta_{s_h^0}(t+h)\mid T=h]\le[R-C]_+.
 \label{eq:bp-retention-cold-value}
\end{equation}
The cached value lies between $R-\varepsilon(L)$ and $R$
by \eqref{eq:bp-retention-warm-value}. Subtracting the two
value bounds gives
\begin{equation}
 \min\{C,R\}-\varepsilon(L)\le g_h\le R,
 \label{eq:bp-retention-gain-bounds}
\end{equation}
where the lower bound uses
$R-\varepsilon(L)-[R-C]_+=\min\{C,R\}-\varepsilon(L)$.

For the short-cache comparison, we can sharpen the upper bound.
For each real route $f$, define its execution saving by
\[
 \delta_{h,f}=\lambda_f(t+h)
       \bigl(c_{s_h^0,f}(t+h)-c_{s_h,f}(t+h)\bigr),
 \qquad \delta_{h,0}=0.
\]
Under the additional assumption, $0\le\delta_{h,f}\le D$.
For the same action, cached and cache-free states have equal
admission rewards, memory charges, and successor values.
Consequently, \eqref{eq:bp-action-value} gives
$A_{s_h,f}(t+h)=A_{s_h^0,f}(t+h)+\delta_{h,f}$, and hence
\[
 \eta_{s_h^0}(t+h)
 \le\max_f\{A_{s_h^0,f}(t+h)+\delta_{h,f}\}
 =\eta_{s_h}(t+h)
 \le\eta_{s_h^0}(t+h)+D.
\]
Taking conditional expectations and using
\eqref{eq:bp-retention-gain-bounds} yields
\begin{equation}
 0\le g_h\le\min\{R,D\}.
 \label{eq:bp-retention-small-gain}
\end{equation}

\paragraph{Step 2: Calculate the value of a waiting deadline.}
Before the next prompt, a deterministic policy observes only
continued nonarrival and the resulting idle age. Along this
no-return history, it has a first eviction date, which defines
a deadline $u\in\{0,\ldots,N\}$. Deadline $u=0$ means
immediate eviction. A positive deadline keeps the cache available
for prompt epochs $t+1,\ldots,t+u$, and evicts it at the
boundary ending $t+u$ if no prompt has arrived. Waiting beyond
the final prompt epoch provides no terminal value.

Write $p_h=\Pp(T=h)$ and omit terms with $p_h=0$.
For $u\ge1$, define
\[
 F(u)=\Pp(T\le u)=\sum_{h=1}^u p_h,
 \qquad
 W(u)=\E[\min(T,u)]=\sum_{h=1}^u\Pp(T\ge h).
\]
Here $F(u)$ is the probability of reuse before eviction, and
$W(u)$ is the expected number of charged waiting boundaries.

Let $U_u$ denote the expected completions minus all resource
charges from the current boundary through $H$, when the controller
uses deadline $u$ and acts optimally after the next prompt.
Immediate eviction has value
\(
 U_0=\sum_{h=1}^N p_h\E[\eta_{s_h^0}(t+h)\mid T=h]
     =v_{s^0}(t+1).
\)
Under a positive deadline, the return state is cached
when $h\le u$ and cache-free when $h>u$. Thus
\begin{align*}
 U_u={}&-\gamma b_s W(u)
       +\sum_{h=1}^u p_h\E[\eta_{s_h}(t+h)\mid T=h]\\
      &+\sum_{h=u+1}^N p_h\E[\eta_{s_h^0}(t+h)\mid T=h].
\end{align*}
Subtracting $U_0$ cancels the terms for $h>u$ and leaves
\begin{equation}
 B_u:=U_u-U_0
 =\sum_{h=1}^u p_hg_h-\gamma b_s W(u).
 \label{eq:bp-retention-deadline}
\end{equation}
Thus, $B_u$ is the gain over immediate eviction with deadline $u$: the expected
benefit of returning before the deadline, less the storage
charge paid while waiting. 

\paragraph{Step 3: Recover the theorem's retention comparison.}
Every deterministic waiting policy corresponds to a deadline.
Fixing a randomized policy's independent random choices also
produces a deadline; averaging over these choices cannot exceed
the best deadline value. The optimal continuation used in
$U_u$ is attainable at every return state by the backward
recursion. Therefore, the best value conditional on KEEP now
is $\max_{1\le u\le N}U_u$, while DROP has value $U_0$.
The boundary recursion expresses these same two values as
\[
 \max_{1\le u\le N}U_u=v_s(t+1)-\gamma b_s,
 \qquad U_0=v_{s^0}(t+1),
\]
where the first equation is because keeping the cache at the current boundary incurs $\gamma b_s$
and leaves the session in state $s$ at the beginning of epoch
$t+1$
Subtracting and recalling
$\Delta_s(t+1)=v_s(t+1)-v_{s^0}(t+1)$ gives
\begin{align}
 \Delta_s(t+1)-\gamma b_s
 &=\max_{1\le u\le N}B_u,\nonumber\\
 w_s(t)-v_{s^0}(t+1)
 &=\max_{0\le u\le N}B_u.
 \label{eq:bp-retention-boundary-bridge}
\end{align}
In particular, a positive $B_u$ for some deadline requires
KEEP; strictly negative $B_u$ for every positive deadline
requires DROP by Theorem~\ref{prop:bp-control}.b.

\paragraph{Step 4: Obtain the three sufficient conditions.}
For any deadline with $F(u)>0$, substitute the lower bound
in \eqref{eq:bp-retention-gain-bounds} into
\eqref{eq:bp-retention-deadline} to obtain
\[
 B_u\ge
 F(u)\left\{\min\{C,R\}-\varepsilon(L)
                    -\gamma b_s\frac{W(u)}{F(u)}\right\}.
\]
Choose $u_*$ minimizing $W(u)/F(u)$ over these deadlines.
Such a minimizer exists because the set is finite and nonempty,
and its ratio equals $\deff$. Hence
\[
 B_{u_*}\ge F(u_*)
       \{\min\{C,R\}-\varepsilon(L)-\gamma b_s\deff\}>0
\]
under \eqref{eq:bp-retention-keep}. Step 3 therefore requires
retention.

For eviction, we must bound the gain of every positive deadline.
If $F(u)>0$, the upper bound $g_h\le R$ gives
\[
 B_u\le F(u)\left\{R-\gamma b_s\frac{W(u)}{F(u)}\right\}
      \le F(u)\{R-\gamma b_s\deff\}<0
\]
under \eqref{eq:bp-retention-long-drop}. For the short-cache
condition, the sharper bound $g_h\le D$ gives instead
\[
 B_u\le F(u)\{D-\gamma b_s\deff\}<0
 \quad\text{if}\quad D<\gamma b_s\deff.
\]
Using $D=\bar\lambda_e(a_1L+a_2L^2)$ and $b_s=\beta L$,
and dividing by $L>0$, the last condition is exactly
\eqref{eq:bp-retention-short-drop}.

Finally, if $F(u)=0$, no return occurs by $u$, so $W(u)=u$
and $B_u=-\gamma b_su<0$. Thus either eviction condition
makes all finitely many positive-deadline gains strictly
negative. Their maximum is negative, and Step 3 requires eviction.
\end{proof}

\subsubsection{Evaluating the deficit and protected lengths}
\label{app:bp-retention-primitive}

Consider the policy that serves the returning prompt at its
owner, continues until completion or $H$, and then drops its
cache and rejects all later prompts. Conditional on $T=h$,
let $q_e(h,L)$ be the probability that this turn completes
by $H$, and let $C_e^{\mathrm w}(h,L)$ be its expected
execution and active-memory charges through $H$. Its value
is $q_e(h,L)-C_e^{\mathrm w}(h,L)$, so evaluating this policy
in \eqref{eq:bp-retention-epsilon} gives
\begin{equation}
 \varepsilon(L)\le
 \max_{\substack{1\le h\le N\\\Pp(T=h)>0}}
       \{R-q_e(h,L)+C_e^{\mathrm w}(h,L)\}.
 \label{eq:bp-retention-one-turn-epsilon}
\end{equation}
Suppose, uniformly over the lengths and dates being compared,
\begin{equation}
 q_e(h,L)\ge q_0>0,\qquad
 C_e^{\mathrm w}(h,L)\le c_0+c_1L,
 \qquad c_0,c_1\ge0,
 \label{eq:bp-retention-warm-bounds}
\end{equation}
where $C_e^{\mathrm w}(h,L)\le c_0+c_1L$ is justified later in this subsection.
Then $q_0\le\min\{1,R\}$ and
$\varepsilon(L)\le R-q_0+c_0+c_1L$. Substitution into
\eqref{eq:bp-retention-keep}, using
$R-\min\{C,R\}=[R-C]_+$, proves the sufficient condition
\begin{equation}
 c_0+(c_1+\beta\gamma\deff)L
 +[R-\underline\lambda(a_1L+a_2L^2)]_+<q_0.
 \label{eq:bp-retention-explicit}
\end{equation}

\paragraph{An explicit interval.}
Hold the return law and these coefficients fixed across
lengths. Suppose $\underline\lambda>0$ and $q_0>c_0$.
Define
\[
 L_c=\frac{-a_1+\sqrt{a_1^2+4a_2R/\underline\lambda}}
                 {2a_2},
 \qquad
 L_+=\frac{q_0-c_0}{c_1+\beta\gamma\deff}.
\]
Here $L_c>0$ solves
$\underline\lambda(a_1L_c+a_2L_c^2)=R$.
If $L_c<L_+$, every admissible length in $[L_c,L_+)$
satisfies \eqref{eq:bp-retention-explicit}: its positive-part
term is zero and its remaining cost is strictly below $q_0$.
This gives the interval of protected lengths stated in the
main text. 

Under the additional short-cache condition and
$\bar\lambda_e>0$, the sufficient eviction interval is
\begin{equation}
 0<L<\frac{\beta\gamma\deff-\bar\lambda_e a_1}
                {\bar\lambda_e a_2},
 \qquad \beta\gamma\deff>\bar\lambda_e a_1.
 \label{eq:bp-retention-short-interval}
\end{equation}
This follows by solving \eqref{eq:bp-retention-short-drop}
for $L$. If $\bar\lambda_e=0$, that test instead implies
eviction at every positive length.

\paragraph{A service-cost example.}
An incremental work model with $m$ new tokens has profiles
\[
 \begin{aligned}
 \text{cold:}\quad&a_0+a_1(L+m)+a_2(L+m)^2,\\
 \text{warm:}\quad&a_0+a_1m+a_2(2Lm+m^2).
 \end{aligned}
\]
Their difference is $a_1L+a_2L^2$: tokenwise work contributes
the linear term and attention within the old prefix contributes
the quadratic term. New tokens still attend to that prefix
when it is cached~\citep{vaswani,hf-cache}.

Suppose the warm turn uses at most $m$ new tokens, its total
execution work is bounded by the warm profile, and its
active footprint is at most $\beta(L+m)$. If
$\bar\lambda_e$ and $\bar\gamma_e$ bound the execution
and active-memory prices over the remaining horizon, then
\begin{align}
 C_e^{\mathrm w}(h,L)
 &\le\bar\lambda_e[a_0+a_1m+a_2(2Lm+m^2)]
       +\bar\gamma_e\beta H(L+m)
   =c_0+c_1L,\label{eq:bp-retention-warm-derived}\\
 c_0&=\bar\lambda_e(a_0+a_1m+a_2m^2)
            +\bar\gamma_e\beta Hm,\nonumber\\
 c_1&=2\bar\lambda_e a_2m+\bar\gamma_e\beta H.\nonumber
\end{align}
The completion bound $q_0$ is evaluated separately from the
service kernels. Both the warm upper bound and cold lower
bound count work actually charged by $H$.

The same policy-evaluation argument allows several future
turns. If a specified policy has expected completions at least
$Q$ and expected charges at most $c_0+c_1L$ at every possible
return date, then
$\varepsilon(L)\le R-Q+c_0+c_1L$, and
\eqref{eq:bp-retention-explicit} holds with $Q$ replacing $q_0$.

\subsubsection{Return distributions and the age comparison}
\label{app:bp-retention-age}

For two residual return times $T_1\le_{\mathrm{st}}T_2$ (meaning that
$\Pp(T_1>u)\le\Pp(T_2>u)$ for every $u$)
over the same remaining horizon, every deadline satisfies
\[
 \E[\min(T_1,u)]\le\E[\min(T_2,u)],\qquad
 \Pp(T_1\le u)\ge\Pp(T_2\le u).
\]
Thus each positive-denominator ratio for $T_2$ is weakly
smaller for $T_1$, and $T_1$ may have additional eligible
deadlines. Minimization gives
$\deff(T_1)\le\deff(T_2)$, with the same infinity convention.

\begin{proof}[Proof of Corollary~\ref{cor:bp-retention-age}]
At an age $a$ with $\Pp(U>a\mid L)>0$, define
\[
 S_a(h)=\Pp(U>a+h\mid U>a,L),\qquad S_a(0)=1.
\]
Let $T_a$ be the residual return time, with returns beyond
the common horizon $N$ recorded as $\infty$. For $h\le N$,
\(
 \Pp(T_a=h)=S_a(h-1)p(a+h-1,L).
\)
The
identity $\min(T_a,u)=\sum_{h=1}^u\1\{T_a\ge h\}$
gives, for $1\le u\le N$,
\begin{align*}
 \E[\min(T_a,u)]&=\sum_{h=1}^u S_a(h-1),\\
 \frac{\Pp(T_a\le u)}{\E[\min(T_a,u)]}
 &=\frac{\sum_{h=1}^u S_a(h-1)p(a+h-1,L)}
         {\sum_{h=1}^u S_a(h-1)}.
\end{align*}
The denominator is positive. If the hazard is nonincreasing,
this weighted average is at most $p(a,L)$, with equality
at $u=1$. Taking reciprocals of the maximum ratio gives
$\deff(a,L)=1/p(a,L)$ when $p(a,L)>0$. If $p(a,L)=0$,
all subsequent hazards vanish, so both sides are infinite.
The expression is nondecreasing in age.

If the hazard is nondecreasing, then
\(
 \Pp(T_a>h)=\prod_{r=0}^{h-1}(1-p(a+r,L)),
 \; h\le N.
\)
Every factor is weakly smaller at an older age, so its
residual return time is stochastically smaller. The ordering
proved above makes $\deff$ nonincreasing in age. Constant
hazard gives $\deff=1/p$, independently of age.
\end{proof}

\subsubsection{Time-varying waiting prices}
\label{app:bp-retention-extensions}

For a deterministic nonnegative owner-price path, replace
$\gamma\deff$ in the three certificates by
\begin{equation}
 \min_{\substack{1\le u\le N\\\Pp(T\le u)>0}}
 \frac{\E\!\left[
       \sum_{h=1}^{\min(T,u)}\gamma_e(t+h-1)\right]}
      {\Pp(T\le u)}.
 \label{eq:bp-retention-priced-wait}
\end{equation}
To verify this, denote the numerator by $W_\gamma(u)$.
The deadline calculation becomes
\[
 B_u=\sum_{h=1}^u p_hg_h-b_sW_\gamma(u).
\]
Using the lower bound on $g_h$ at a minimizing deadline
proves the retention certificate. Using either upper bound
at every deadline with $F(u)>0$ proves the corresponding
eviction certificate. A deadline with $F(u)=0$ has gain
$-b_sW_\gamma(u)\le0$. Hence eviction is optimal under
either strict eviction inequality; it is the unique initial
action if $\gamma_e(t)>0$, because every positive deadline
then incurs a strictly positive charge.

If warm-service bounds are available only for returns in
the first $m\le N$ epochs, the retention test can restrict
both the deficit's maximum and the deadline minimum to this
range. Its proof chooses a deadline $u\le m$ and therefore
uses only these bounds. Eviction tests continue to use all
possible return dates and deadlines.

\section{Proofs for the buffer-reserve policy}
\label{app:br-proofs}

\paragraph{Notation and the emergency rule.}
Fluid populations and flows in this appendix are divided by $n$.
Write
\[
 \delta_n=r_n/n=b\sqrt{\log n/n},
 \qquad
 \mathcal W_t(m;\mathbf c)=N_A(m)+\mathcal V_t(m;\mathbf c),
\]
where $\mathbf c$ collects execution and memory capacities over the
remaining horizon. The nominal profile is $\mathbf c_t$; the profile
used by the buffered LP at $t$ is $\mathbf c_t^{-,n;t}$, with entries
$B_e-(\tau-t+1)\delta_n$ and $G_e-(\tau-t+1)\delta_n$ for $\tau\ge t$.
Thus $\mathcal W_t(m)$ in the main text means
$\mathcal W_t(m;\mathbf c_t)$. An implemented quota is
$\lfloor ny_{s,e}\rfloor$ or $\lfloor nz_{s,\KEEP}\rfloor$ in these units.

For the analysis, the \emph{potential ordinary epoch} follows the
quota decisions and the ordinary service laws without intervention.
The theoretical emergency mechanism in Section~3 enforces the actual
capacity constraints causally. We use its following property: it
leaves the epoch unchanged whenever the potential ordinary execution
loads and post-retention boundary-memory loads fit. Otherwise it
can stop admissions, interrupt unfinished work, and release memory
before the corresponding constraint is violated. An infeasible
reoptimization also triggers the emergency rule. Previously completed
turns keep their rewards, and interrupted turns earn none. This
specifies the safety mechanism used in the theorem; future potential
loads are used only in its probability analysis.

Let $\bar\nu=\max_\tau\sum_s\nu_s(\tau)$ and
$K_\tau=K_0+\tau(\bar\nu+1)$. The estimates below are uniform over
histories with $\|M^{(n)}(\tau)\|_1\le K_\tau n$. On the good
successor events used below, the next population satisfies the
corresponding bound with $K_{\tau+1}$. Constants may depend on
$K_0,H$, the fixed primitives, and $b$, but not on $n$ or the history.

\subsection{Two properties of the fluid LP}
We first establish a uniform fluid-value regularity.
\begin{lemma}
\label{lem:br-lipschitz}
There is a finite constant $L$, uniform over starting epochs, such
that any two feasible parameter pairs satisfy
\[
 |\mathcal W_t(m;\mathbf c)-\mathcal W_t(m';\mathbf c')|
 \le L\bigl(\|m-m'\|_1+\|\mathbf c-\mathbf c'\|_1\bigr).
\]
The same conclusion holds for $\mathcal V_t$. Both values are
nondecreasing in each capacity coordinate.
\end{lemma}

\begin{proof}
Sub-stochastic transitions and flow conservation give
\[
 \|m(\tau+1)\|_1
 \le\|m(\tau)\|_1+\sum_s\nu_s(\tau).
\]
This bounds every population variable over the finite horizon.
Arrival and retention conservation then bound every action variable.
Each feasible instance consequently has a finite attained optimum.

Add resource slacks, remove linearly dependent equality rows, and
write the remaining LP as
\[
 \max\{f^\top x:Ax=d+Dm+F\mathbf c,\ x\ge0\},
\]
where the matrices are fixed and $A$ has full row rank. 
The dual feasible set $\{v:A^\top v\ge f\}$ is nonempty by LP
strong duality. This polyhedron has no line: a
line direction would satisfy $A^\top v=0$, and full row rank forces
$v=0$. It therefore has finitely many vertices
$v^1,\ldots,v^J$, and a finite attained dual optimum occurs at a
vertex. For every feasible parameter pair,
\[
 \mathcal V_t(m;\mathbf c)
 =\min_{1\le j\le J}(v^j)^\top(d+Dm+F\mathbf c).
\]
Let $L_0$ be the largest infinity norm of the slope vectors $((v^j)^\top D,(v^j)^\top F)$ of these
affine functions. If $j'$ minimizes at $(m',\mathbf c')$, evaluation of the same affine
function at $(m,\mathbf c)$ gives
\[
 \mathcal V_t(m;\mathbf c)-\mathcal V_t(m';\mathbf c')
 \le L_0(\|m-m'\|_1+\|\mathbf c-\mathbf c'\|_1).
\] 
Interchanging the
pairs proves the absolute bound. Adding $N_A(m)$ increases the
constant by at most one. Enlarging capacities enlarges the feasible
set, proving monotonicity. Taking the largest constant over the
finitely many starting epochs gives $L$.
\end{proof}

Next we establish a recourse for a perturbed initial population.
\begin{lemma}
\label{lem:br-recourse}
If the remaining-horizon fluid LP is feasible from $m^1$ under
$\mathbf c$, then it is feasible from any $m^2\ge0$ under
\[
 \mathbf c+L_{\mathrm{rec}}\|m^2-m^1\|_1\mathbf1,
 \qquad
 L_{\mathrm{rec}}
 =\max_{\tau,e,s\in\cSA}
   \{\ell_{s,e}(\tau),\kappa_{s,e}^{\mathrm{act}}(\tau)\}.
\]
The maximum covers the remaining horizon and is zero if there are
no active states.
\end{lemma}

\begin{proof}
Start with a feasible solution from $m^1$. At each idle state and
epoch, extract its routing fractions
\[
 \alpha_{s,e}(\tau)
 =\frac{y^1_{s,e}(\tau)}{p_s(\tau)m_s^1(\tau)+\nu_s(\tau)}
\]
when the denominator is positive; otherwise assign fraction one to
rejection. Similarly, extract the retention fractions
\[
 \beta_{s,h}(\tau)
 =\frac{z^1_{s,h}(\tau)}{\widetilde m_s^1(\tau+1)},
 \qquad h\in\{\KEEP,\DROP\},
\]
using DROP with fraction one when the denominator is zero.
Each set of fractions is nonnegative and sums to one.

Let $m_s^\wedge=\min\{m_s^1,m_s^2\}$. Starting from $m^\wedge$,
apply these fractions to the new populations, including all original
exogenous sources. Specifically, set
\[
 \widehat y_{s,e}(\tau)
 =\alpha_{s,e}(\tau)
       [p_s(\tau)\widehat m_s(\tau)+\nu_s(\tau)],
 \qquad
 \widehat z_{s,h}(\tau)
 =\beta_{s,h}(\tau)\widehat{\widetilde m}_s(\tau+1),
\]
and obtain the populations from the pre- and post-retention equations.
All flow equations hold by construction. If
$\widehat m(\tau)\le m^1(\tau)$ componentwise, the same fractions
and exogenous sources imply $\widehat y(\tau)\le y^1(\tau)$.
Nonnegative transition coefficients then imply
$\widehat{\widetilde m}(\tau+1)\le\widetilde m^1(\tau+1)$,
followed by $\widehat z(\tau)\le z^1(\tau)$ and
$\widehat m(\tau+1)\le m^1(\tau+1)$. These inequalities hold
initially and hence throughout the horizon. Every resource load is
therefore no greater than in the original solution.

Add a separate flow of mass $d=m^2-m^\wedge$, with no exogenous
sources. Its idle sessions reject every prompt and drop every cache.
Its initially active turns continue under $Q^A$; after completion,
they use the same reject-and-drop rule. No new active turn is created,
and sub-stochastic transitions cannot increase the total mass.
Consequently, the added active mass never exceeds
\[
 \|d\|_1\le\|m^2-m^1\|_1.
\]
Its load in any execution or memory row is at most
$L_{\mathrm{rec}}\|m^2-m^1\|_1$. Idle caches contribute no load
because memory is measured after boundary eviction. Superposing this
flow on the solution from $m^\wedge$ gives initial population $m^2$,
retains exactly the prescribed exogenous sources, and satisfies all
flow equations. The stated capacity increase covers every added load.
\end{proof}

\subsection{Quota errors and the successor population}

Fix a feasible buffered LP at epoch $t$ and condition on its beginning
history $\mathcal F_t$, with $M=nm$. Suppress $t$ in the normalized
current flows $(y,z)$, and write $\widetilde m^+,m^+$ for their
pre- and post-retention successors. In the potential ordinary epoch,
let $Y_{s,e}$ count action-$e$ arrivals from state $s$,
$\widetilde M_s$ count pre-retention sessions, and $M^+$ be the
post-retention population. The boundary rule is
\begin{equation}
 Z_{s,\KEEP}
 =\min\{\lfloor nz_{s,\KEEP}\rfloor,\widetilde M_s\},
 \qquad
 Z_{s,\DROP}=\widetilde M_s-Z_{s,\KEEP},
 \qquad s\in\cSI,
 \label{eq:br-boundary-quota}
\end{equation}
with zero KEEP flow in cache-free states.

We establish some concentration results for quota and successor.
\begin{lemma}
\label{lem:br-concentration}
There are constants $C,c>0$ such that
\begin{align}
 \E\!\left[
   \sum_{s\in\cSI,e\in\cEzero}|Y_{s,e}-ny_{s,e}|
       +\|M^+-nm^+\|_1\,\middle|\,\mathcal F_t\right]
 &\le C\sqrt n,\label{eq:br-mean-error}\\
 \Pp\!\left(
   \|M^+-nm^+\|_1>C+u\mid\mathcal F_t\right)
 &\le C\exp\!\left\{-\frac{u^2}{c(n+u)}\right\},
 \quad u\ge0.
 \label{eq:br-state-tail}
\end{align}
The route-count error has the  tail
$C\exp\{-u^2/(cn)\}$ above a constant offset, i.e.,
\[
 \Pp\!\left(
   \sum_{s\in\cSI,e\in\cEzero}|Y_{s,e}-ny_{s,e}|
   >2E|\cSI|+u
   \,\middle|\,\mathcal F_t
 \right)
 \le C\exp\!\left\{-\frac{u^2}{cn}\right\},
 \qquad u\ge0.
\]
\end{lemma}

\begin{proof}
Let $A_s$ be endogenous arrivals, $D_s=A_s+\Xi_s^{(n)}$, and
$q_s=p_sm_s+\nu_s$. Define
\[
 U_A=\sum_{s\in\cSI}|A_s-np_sm_s|,
 \qquad U_D=\sum_{s\in\cSI}|D_s-nq_s|.
\]
Conditional Bernoulli independence, $M_s\le K_tn$, and Hoeffding's
inequality bound each endogenous deviation by a constant times
$\exp\{-u^2/(cn)\}$. Assumption~\ref{ass:br-scaling} gives the
same bound for each exogenous deviation. If a sum of $d$ absolute
deviations exceeds $u$, one exceeds $u/d$; a union bound over the
fixed number of states therefore gives
\begin{equation}
 \Pp(U_A+U_D>u\mid\mathcal F_t)
 \le C\exp\{-u^2/(cn)\}.
 \label{eq:br-arrival-tail}
\end{equation}

For each state, let
$Q_s=\sum_{e\in\cE}\lfloor ny_{s,e}\rfloor$ and
$a_s=\sum_{e\in\cE}y_{s,e}\le q_s$. Sequential quota filling gives
\[
 \sum_{e\in\cE}Y_{s,e}=\min\{D_s,Q_s\},
 \qquad Y_{s,e}\le ny_{s,e}\quad(e\in\cE).
\]
Since $na_s-E\le Q_s\le na_s$,
\[
 \sum_{e\in\cE}|Y_{s,e}-ny_{s,e}|
 =na_s-\min\{D_s,Q_s\}
 \le E+|D_s-nq_s|.
\]
Arrival conservation also gives
$|Y_{s,0}-ny_{s,0}|\le E+2|D_s-nq_s|$. Thus
\begin{equation}
 U_Y:=\sum_{s\in\cSI,e\in\cEzero}|Y_{s,e}-ny_{s,e}|
 \le2E|\cSI|+3U_D.
 \label{eq:br-route-total-bound}
\end{equation}

Condition further on the arrival list, writing
$\mathcal G_t=\mathcal F_t\vee\sigma(\text{arrival list})$.
The quota routes are then fixed, and
\[
 \mu_r:=\E[\widetilde M_r\mid\mathcal G_t]
 =\sum_{s\in\cSA}M_sQ_t^A(s,r)
  +\sum_{s\in\cSI}(M_s-A_s)Q_t^I(s,r)
  +\sum_{s\in\cSI,e\in\cEzero}Y_{s,e}P_t^e(s,r).
\]
At most $N_t=\|M\|_1+\sum_s\Xi_s^{(n)}$ independent categorical
outcomes contribute to $\widetilde M$. The new-session count need
not have a deterministic bound proportional to $n$. Instead, choose
$a>K_t+\sum_s\nu_s(t)+1$. The exogenous tail bound implies
\[
 \Pp(N_t>an+u\mid\mathcal F_t)
 \le C\exp\{-(n+u)^2/(cn)\}.
\]
On any arrival list with $N_t\le an+u$, Hoeffding's inequality
for each coordinate, followed by a union bound, yields
\[
 \Pp(\|\widetilde M-\mu\|_1>u\mid\mathcal G_t)
 \le C\exp\{-u^2/[c(n+u)]\}.
\]
Splitting according to $N_t\le an+u$ and taking expectations proves
\[
 \Pp(\|\widetilde M-\mu\|_1>u\mid\mathcal F_t)
 \le C\exp\{-u^2/[c(n+u)]\}.
\]
Sub-stochastic kernel rows have one-norm at most one, so subtracting
the fluid transition equation gives
\[
 \|\mu-n\widetilde m^+\|_1\le U_A+U_Y.
\]
Combining the last two bounds with
\eqref{eq:br-arrival-tail}--\eqref{eq:br-route-total-bound} shows that
$U_{\mathrm{pre}}:=\|\widetilde M-n\widetilde m^+\|_1$ has the
stated tail above a fixed constant offset. This also proves the claimed concentration bound for routing
counts. Indeed, for every $u\ge0$,
\begin{align*}
 \Pp\!\left(U_Y>2E|\cSI|+u\mid\mathcal F_t\right)
 &\le \Pp\!\left(U_D>u/3\mid\mathcal F_t\right)\\
 &\le \Pp\!\left(U_A+U_D>u/3\mid\mathcal F_t\right)\\
 &\le C\exp\!\left\{-\frac{u^2}{9cn}\right\}.
\end{align*}
The first inequality follows from
\eqref{eq:br-route-total-bound}, the second uses $U_A\ge0$,
and the third applies \eqref{eq:br-arrival-tail} with
threshold $u/3$. Absorbing the factor $9$ into the constant
$c$ gives the stated routing-count bound.

Because $z_{s,\KEEP}\le\widetilde m_s^+$, the minimum in
\eqref{eq:br-boundary-quota} implies
\begin{align*}
 |Z_{s,\KEEP}-nz_{s,\KEEP}|
 &\le1+|\widetilde M_s-n\widetilde m_s^+|,\\
 |Z_{s,\DROP}-nz_{s,\DROP}|
 &\le1+2|\widetilde M_s-n\widetilde m_s^+|.
\end{align*}
The second inequality follows by subtracting KEEP from retention
conservation. Active mass passes unchanged through retention, and
DROP sends each idle column to one state $\Ev(s)$. Consequently,
\[
 \|M^+-nm^+\|_1\le2|\cSI|+3U_{\mathrm{pre}},
\]
which proves \eqref{eq:br-state-tail}.

For a nonnegative variable with this tail, its mean is at most a
constant plus
\[
 C\int_0^n e^{-u^2/(2cn)}\,du
 +C\int_n^\infty e^{-u/(2c)}\,du
 =O(\sqrt n).
\]
Here $n+u\le2n$ on the first interval and $n+u\le2u$ on the
second. Applying this bound to the successor and route errors proves
\eqref{eq:br-mean-error}.
\end{proof}

\subsection{Current safety and feasibility of the next plan}
The following lemma shows that a sufficiently large capacity reserve absorbs single-epoch stochastic fluctuations without emergency intervention and preserves feasibility of the buffered linear program across the horizon with high probability.
\begin{lemma}
\label{lem:br-rare}
For any $q>0$, choose $b$ (from (\ref{eq:br-buffer-size})) sufficiently large. At a history with a
feasible buffered LP and $\|M^{(n)}(t)\|_1\le K_tn$, there is an
event of conditional probability at least $1-Cn^{-q}$ on which the
current epoch needs no intervention and the next buffered LP is
feasible. On this event, $\|M^{(n)}(t+1)\|_1\le K_{t+1}n$.
Starting from a feasible buffered LP and $\|M^{(n)}(t)\|_1\le K_0n$,
the probability of an emergency during the remaining horizon is at
most $Cn^{-q}$.
\end{lemma}

\begin{proof}
We first show that the reserve covers current resource
fluctuations, then show that the next buffered LP remains
feasible, and finally combine these bounds over the horizon. All current-epoch coefficients suppress their argument $t$.

\paragraph{Step 1: Current resource use stays within capacity.}
Fix the beginning-of-epoch history $\mathcal F_t$, and let
\(
 \mathcal G_t
 =\mathcal F_t\vee\sigma(\text{arrival list})
\)
include the epoch's arrival list. Conditional on $\mathcal G_t$,
the quota rule determines which prompts are admitted and their
assigned endpoints. 

Consider the epoch under ordinary service, without intervention.
At endpoint $e$, its execution load and unfinished-turn memory are $X_e
 =\sum_{i\in\mathcal A_t}L_{t,i,e}
   +\sum_{\substack{1\le j\le J_t:\\e_{t,j}=e}}C_{t,j,e},A_e
 =\sum_{i\in\mathcal A_t}K_{t,i,e}^{\mathrm{act}}
   +\sum_{\substack{1\le j\le J_t:\\e_{t,j}=e}}
                      K_{t,j,e}^{\mathrm{adm}}.$
Thus, each summand is the resource use of one initially active
session or one newly admitted prompt. The memory terms are zero
for turns that complete before the boundary. Let $I_e$ be the
memory occupied by idle caches retained at that boundary.
Using the conditional resource means gives
\begin{align*}
 \mu_e^X:=\E[X_e\mid\mathcal G_t]
 &=\sum_{s\in\cSA}\ell_{s,e}M_s
   +\sum_{s\in\cSI}c_{s,e}Y_{s,e},\\
 \mu_e^A:=\E[A_e\mid\mathcal G_t]
 &=\sum_{s\in\cSA}\kappa_{s,e}^{\mathrm{act}}M_s
   +\sum_{s\in\cSI}\kappa_{s,e}^{\mathrm{adm}}Y_{s,e}.
\end{align*}
For each real endpoint, the quota rule ensures
$Y_{s,e}\le ny_{s,e}$. Define the fluid retention footprint
\(
 \overline I_e
 =n\sum_{\substack{s\in\cSI:\\o(s)=e}}b_sz_{s,\KEEP}.
\)
The buffered LP reserves $r_n$ units in each current resource
constraint. Since resource coefficients are nonnegative,
\begin{equation}
 \mu_e^X\le nB_e-r_n,
 \qquad
 \mu_e^A+\overline I_e\le nG_e-r_n.
 \label{eq:br-resource-slack}
\end{equation}
Moreover, the retention quotas imply $I_e\le\overline I_e$
on every realization.

We next bound deviations from these conditional means.
The total number of sessions receiving service is
\(
 N_t^{\mathrm{srv}}
 :=|\mathcal A_t|
   +\sum_{s\in\cSI,e\in\cE}Y_{s,e}.
\)
This number is fixed conditional on $\mathcal G_t$. Arrival
conservation in the fluid LP and the routing quotas give
\begin{align*}
 N_t^{\mathrm{srv}}
 &\le \|M\|_1+n\sum_{s\in\cSI,e\in\cE}y_{s,e}
 \le \|M\|_1+n\sum_{s\in\cSI}(p_sm_s+\nu_s)\\
 &\le \left(2K_t+\sum_{s\in\cSI}\nu_s\right)n
 \le C_Nn,
\end{align*}
where $C_N$ is independent of $n$. 

By the bounded-resource condition, each session's execution
requirement and boundary-memory footprint lie in $[0,K]$
for a fixed $K$. Conditional on $\mathcal G_t$, these resource
outcomes are independent across sessions. Consequently, either
$W=X_e$ or $W=A_e$ is a sum of at most $N_t^{\mathrm{srv}}$
independent bounded random variables. Hoeffding's inequality
therefore gives
\[
 \Pp\!\left(
 W-\E[W\mid\mathcal G_t]>u\mid\mathcal G_t
 \right)
 \le
 \exp\!\left\{-\frac{2u^2}{N_t^{\mathrm{srv}}K^2}\right\}
 \le
 \exp\!\left\{-\frac{u^2}{C_Rn}\right\},
 \qquad u>0,
\]
for a constant $C_R$. 

The reserved slack converts these deviation bounds into
capacity guarantees. In particular,
\(
 \{X_e>nB_e\}\subseteq\{X_e-\mu_e^X>r_n\}.
\)
For memory, the pathwise retention bound gives
\[
 A_e+I_e
 \le A_e+\overline I_e
 \le nG_e-r_n+(A_e-\mu_e^A),
\]
so
\(
 \{A_e+I_e>nG_e\}
 \subseteq\{A_e-\mu_e^A>r_n\}.
\)
This argument allows retention to depend on the observed service
outcomes: only the upper bound $I_e\le\overline I_e$ is needed.

Apply the concentration inequality with $u=r_n$, average over
the arrival list, and take a union bound over the execution and
memory constraints at all endpoints. Since
$r_n=b\sqrt{n\log n}$,
\begin{equation}
 \Pp\!\left(
   \text{some ordinary resource load exceeds capacity}
   \mid\mathcal F_t
 \right)
 \le 2E\exp\!\left\{-\frac{r_n^2}{C_Rn}\right\}
 =2E n^{-b^2/C_R}.
 \label{eq:br-resource-tail}
\end{equation}
When all these loads fit, the safety mechanism makes no
intervention.

\paragraph{Step 2: The next buffered LP remains feasible.}
Let $m^+$ be the fluid successor and $M^+$ the successor
under the ordinary stochastic epoch. Set
\[
 K_{\mathrm{rec}}=\max\{1,L_{\mathrm{rec}}\},
 \qquad
 d_n=\frac{\delta_n}{2K_{\mathrm{rec}}},
 \qquad \delta_n=\frac{r_n}{n}.
\]
For sufficiently large $n$, the constant offset in
\eqref{eq:br-state-tail} is at most $nd_n/2$, and
$nd_n/2\le n$. Therefore, applying that inequality with
$u=nd_n/2$ yields
\begin{align}
 \Pp\!\left(
   \|M^+/n-m^+\|_1>d_n\mid\mathcal F_t
 \right)
 &\le C\exp\!\left\{-\frac{nd_n^2}{8c}\right\}
 =C n^{-b^2/(32cK_{\mathrm{rec}}^2)}.
 \label{eq:br-layer-tail}
\end{align}

The future portion of the current fluid solution is feasible
from $m^+$. At each future epoch $\tau\ge t+1$, its execution
capacity is
\(
 B_e-(\tau-t+1)\delta_n,
\)
whereas the LP reoptimized at $t+1$ uses
\(
 B_e-(\tau-t)\delta_n.
\)
The latter therefore has an additional $\delta_n$ in every
execution constraint, and the same increase occurs in every
memory constraint.

On the event $\|M^+/n-m^+\|_1\le d_n$,
Lemma~\ref{lem:br-recourse} constructs a feasible solution
from the realized population using at most
\(
 L_{\mathrm{rec}}\|M^+/n-m^+\|_1
 \le L_{\mathrm{rec}}d_n
 \le\delta_n/2
\)
additional capacity in each constraint. The available increase
$\delta_n$ covers this adjustment, so the next buffered LP
is feasible.

The population bound also continues to hold. Fluid mass can
increase only through new sessions, so
$\|m^+\|_1\le K_t+\bar\nu$. For sufficiently large $n$,
$d_n\le1$, and hence
\(
 \|M^+/n\|_1
 \le\|m^+\|_1+d_n
 \le K_t+\bar\nu+1
 =K_{t+1}.
\)
At $t=H-1$, no further LP is needed; the terminal memory
constraint is already covered by Step 1.

\paragraph{Step 3: Combine the bounds over the horizon.}
Choose $b$ large enough that
\[
 \frac{b^2}{C_R}>q,
 \qquad
 \frac{b^2}{32cK_{\mathrm{rec}}^2}>q.
\]
Call an epoch successful when its ordinary resource loads fit
and its successor error is at most $d_n$. Steps 1 and 2 show
that, conditional on a feasible buffered LP and the population
bound, an epoch is unsuccessful with probability at most
$Cn^{-q}$. A successful epoch requires no intervention and
preserves both conditions for the next epoch.

Let $\mathsf A_\tau$ be the event that all epochs from $t$
through $\tau-1$ are successful. The events identifying the
first unsuccessful epoch are disjoint, so
\begin{align*}
 \Pp(\text{some epoch is unsuccessful})
 &=\sum_{\tau=t}^{H-1}
   \E\!\left[
     \1_{\mathsf A_\tau}
     \Pp(\text{epoch $\tau$ is unsuccessful}\mid\mathcal F_\tau)
   \right]
 \le HCn^{-q}.
\end{align*}
An emergency can occur only if some epoch is unsuccessful.
Since $H$ is fixed, absorbing it into the constant proves
the whole-horizon bound. Steps 1 and 2 also establish the
one-epoch assertion.
\end{proof}

\subsection{The fluid upper bound}

\begin{lemma}
\label{lem:br-upper}
Every feasible nonanticipative policy without service interruptions
has expected remaining completions at most $n\mathcal W_t(m;\mathbf c_t)$
when its initial population is $nm$.
\end{lemma}

\begin{proof}
Fix a feasible nonanticipative policy without service interruptions.
We show that its expected populations and decisions form a feasible
fluid solution with the same expected completion reward, after
division by $n$.

All expectations below condition on the initial population
$M^{(n)}(t)=nm$. Define
\[
 \begin{aligned}
 \bar m_s(\tau)
   &=\E[M_s^{(n)}(\tau)]/n,
 &\overline{\widetilde m}_s(\tau+1)
   &=\E[\widetilde M_s^{(n)}(\tau+1)]/n,\\
 \bar y_{s,e}(\tau)
   &=\E[Y_{s,e}^{(n)}(\tau)]/n,
 &\bar z_{s,h}(\tau)
   &=\E[Z_{s,h}^{(n)}(\tau)]/n.
 \end{aligned}
\]
These expectations are finite. The population is bounded by
the initial population plus all new sessions entering so far,
and each session submits at most one prompt per epoch.

\paragraph{Step 1: Expected flows satisfy the fluid equations.}
Every arriving prompt receives exactly one action, including
rejection. Therefore,
\(
 \sum_{e\in\cEzero}Y_{s,e}^{(n)}(\tau)
 =A_s^{(n)}(\tau)+\Xi_s^{(n)}(\tau).
\)
Taking expectations, using the conditional arrival means,
and dividing by $n$ gives
\(
 \sum_{e\in\cEzero}\bar y_{s,e}(\tau)
 =p_s(\tau)\bar m_s(\tau)+\nu_s(\tau).
\)

Next consider the pre-retention population. Fix an epoch $\tau$,
and define
\[
 I_{j,s,e}
 =\1\{\text{arrival $j$ occurs in state $s$ and receives action $e$}\}.
\]
If fewer than $j$ arrivals occur, this indicator and its
next-state contribution are zero. Let $\widetilde s_j$ denote
the pre-retention state of the session submitting arrival $j$,
and let $\mathcal H_{\tau,j}$ contain the information available
at its routing decision, including the selected action.

At that decision, $I_{j,s,e}$ is already known. Conditional on
$I_{j,s,e}=1$, the probability of next state $r$ is
$P_\tau^e(s,r)$. Hence
\(
 \E\!\left[
   I_{j,s,e}\1\{\widetilde s_j=r\}
   \mid\mathcal H_{\tau,j}
 \right]
 =I_{j,s,e}P_\tau^e(s,r).
\)
Taking expectations and summing over arrivals yields
\begin{align*}
 \E\!\left[\sum_{j\ge1}
   I_{j,s,e}\1\{\widetilde s_j=r\}\right]
 &=\sum_{j\ge1}P_\tau^e(s,r)\E[I_{j,s,e}]
 =P_\tau^e(s,r)\E[Y_{s,e}^{(n)}(\tau)].
\end{align*}
The sums have nonnegative terms and a finite expected number
of arrivals, so expectation and summation can be interchanged.

Initially active sessions contribute
\(
 Q_\tau^A(s,r)\E[M_s^{(n)}(\tau)]
\)
in expectation to state $r$. For idle sessions that submit no
prompt, the corresponding contribution is
\(
 Q_\tau^I(s,r)
 \E[M_s^{(n)}(\tau)-A_s^{(n)}(\tau)]
 =(1-p_s(\tau))Q_\tau^I(s,r)\E[M_s^{(n)}(\tau)].
\)
Adding these contributions and dividing by $n$ gives
\begin{align*}
 \overline{\widetilde m}_r(\tau+1)
 ={}&\sum_{s\in\cSA}\bar m_s(\tau)Q_\tau^A(s,r)
 +\sum_{s\in\cSI}(1-p_s(\tau))
                  \bar m_s(\tau)Q_\tau^I(s,r)
 +\sum_{s\in\cSI,e\in\cEzero}
                  \bar y_{s,e}(\tau)P_\tau^e(s,r).
\end{align*}
This is the fluid pre-retention equation.

At the boundary, every idle session is assigned either KEEP
or DROP. Taking expectations in the realized retention and
state-update equations gives
\begin{align*}
 \bar z_{s,\KEEP}(\tau)+\bar z_{s,\DROP}(\tau)
 &=\overline{\widetilde m}_s(\tau+1),
 &&s\in\cSI,\\
 \bar m_r(\tau+1)
 &=\bar z_{r,\KEEP}(\tau)
   +\sum_{s:\Ev(s)=r}\bar z_{s,\DROP}(\tau),
 &&r\in\cSI,\\
 \bar m_r(\tau+1)
 &=\overline{\widetilde m}_r(\tau+1),
 &&r\in\cSA.
\end{align*}
These identities hold even when retention uses the observed
service outcomes, because the realized equations hold before
taking expectations. 

\paragraph{Step 2: Expected flows satisfy the capacity constraints.}
For a real endpoint $e$, conditioning at the routing decision
also gives
\begin{align*}
 \E[I_{j,s,e}C_{\tau,j,e}]
 &=\E\!\left[
   I_{j,s,e}\E[C_{\tau,j,e}\mid\mathcal H_{\tau,j}]
   \right]
 =c_{s,e}(\tau)\E[I_{j,s,e}].
\end{align*}
Thus, summing over arrivals, the expected execution requirement
of state-$s$ prompts routed to $e$ is
$c_{s,e}(\tau)\E[Y_{s,e}^{(n)}(\tau)]$.
Replacing execution requirements by admission-memory footprints
gives the same equality with
$\kappa_{s,e}^{\mathrm{adm}}(\tau)$.

For initially active sessions, conditioning on the
beginning-of-epoch history gives expected execution and memory
loads
\[
 \sum_{s\in\cSA}\ell_{s,e}(\tau)
                         \E[M_s^{(n)}(\tau)],
 \qquad
 \sum_{s\in\cSA}\kappa_{s,e}^{\mathrm{act}}(\tau)
                         \E[M_s^{(n)}(\tau)].
\]
Finally, the realized retained idle footprint is exactly
$\sum_{s\in\cSI:o(s)=e}b_sZ_{s,\KEEP}^{(n)}(\tau)$.

Taking expectations in the almost-sure capacity constraints
and dividing by $n$ therefore yields
\begin{align*}
 \sum_{s\in\cSA}\ell_{s,e}(\tau)\bar m_s(\tau)
 +\sum_{s\in\cSI}c_{s,e}(\tau)\bar y_{s,e}(\tau)
 &\le B_e,\\
 \sum_{s\in\cSA}\kappa_{s,e}^{\mathrm{act}}(\tau)\bar m_s(\tau)
 +\sum_{s\in\cSI}\kappa_{s,e}^{\mathrm{adm}}(\tau)\bar y_{s,e}(\tau)
 +\sum_{s\in\cSI:o(s)=e}b_s\bar z_{s,\KEEP}(\tau)
 &\le G_e.
\end{align*}
Together with Step 1 and $\bar m(t)=m$, these inequalities
show that the expected flows are feasible for the nominal
fluid LP.

\paragraph{Step 3: The fluid objective equals expected completions.}
Without interruptions, each initially active or newly admitted
turn either completes by $H$ or remains active at $H$.
Consequently, on every realization,
\[
 \text{completed turns}
 =N_A(M^{(n)}(t))
  +\sum_{\tau=t}^{H-1}
       \sum_{s\in\cSI,e\in\cE}Y_{s,e}^{(n)}(\tau)
  -N_A(M^{(n)}(H)),
\]
where $N_A(m)=\sum_{s\in\cSA}m_s$.
Taking expectations gives
\[
 \E[\text{completed turns}]
 =n\left[
   N_A(m)
   +\sum_{\tau=t}^{H-1}
       \sum_{s\in\cSI,e\in\cE}\bar y_{s,e}(\tau)
   -N_A(\bar m(H))
 \right].
\]
The expression in brackets is the completion objective of the
feasible fluid solution constructed above. It cannot exceed
the optimal fluid value $\mathcal W_t(m;\mathbf c_t)$,
which proves the result.
\end{proof}

\subsection{Value tracking and the performance guarantee}
Let us track the buffered completion value.
\begin{lemma}
\label{lem:br-tracking}
Choose the reserve in Lemma~\ref{lem:br-rare} for $q>2$. Uniformly
on the bounded population domains specified above, every feasible
buffered initial state satisfies
\[
 V_{t,n}^{\mathrm{BR}}(nm)
 \ge n\mathcal W_t(m;\mathbf c_t^{-,n;t})
      -A_t\sqrt n-D_tn^{1-q},
\]
for constants $A_t,D_t$ independent of $n$ and the observed history.
\end{lemma}

\begin{proof}
We use backward induction. At $H$ there are no further completions,
and $\mathcal W_H(m)=N_A(m)-N_A(m)=0$, so take $A_H=D_H=0$.

At epoch $t<H$, let $m^+$ be the successor of the chosen optimal
buffered solution. Write $\mathbf c_{\mathrm{old}}$ for its capacity
tail and $\mathbf c_{\mathrm{new}}$ for the next reoptimization's
profile, which is larger by $\delta_n$ in each coordinate. The
current solution's tail is optimal from $m^+$ under
$\mathbf c_{\mathrm{old}}$, since a better feasible tail would
improve the full solution. With
\(
 g_t=N_A(m)+\sum_{s\in\cSI,e\in\cE}y_{s,e}-N_A(m^+),
\)
cancelling the intermediate active count gives
\begin{equation}
 \mathcal W_t(m;\mathbf c_t^{-,n;t})
 =g_t+\mathcal W_{t+1}(m^+;\mathbf c_{\mathrm{old}}).
 \label{eq:br-fluid-telescope}
\end{equation}

In the potential ordinary epoch, its completion count is
\(
 C_t^{\mathrm{pot}}
 =N_A(nm)+\sum_{s\in\cSI,e\in\cE}Y_{s,e}-N_A(M^+).
\)
The quota and successor bounds in Lemma \ref{lem:br-concentration} give
\begin{equation}
 \E[|C_t^{\mathrm{pot}}-ng_t|\mid\mathcal F_t]\le C\sqrt n,
 \qquad
 \E[\|M^+/n-m^+\|_1\mid\mathcal F_t]\le C/\sqrt n.
 \label{eq:br-tracking-errors}
\end{equation}

Let $E_t$ be the good event from Lemma~\ref{lem:br-rare}. On this
event, the actual epoch equals the potential one and the next LP
is feasible. The successor also lies in the next bounded population
domain. Lipschitz continuity in Lemma \ref{lem:br-lipschitz} and capacity monotonicity imply
\begin{align*}
 \mathcal W_{t+1}(M^+/n;\mathbf c_{\mathrm{new}})
 &\ge\mathcal W_{t+1}(m^+;\mathbf c_{\mathrm{new}})
      -L\|M^+/n-m^+\|_1\\
 &\ge\mathcal W_{t+1}(m^+;\mathbf c_{\mathrm{old}})
      -L\|M^+/n-m^+\|_1.
\end{align*}
Both parameter pairs in the first line are feasible; the recourse
lemma establishes this for the realized successor. At $t=H-1$,
the continuation values are zero.

Future completion rewards are nonnegative even on $E_t^c$. Applying
the induction hypothesis on $E_t$ gives the following bound, with
the integrand set to zero outside $E_t$:
\begin{align*}
 V_{t,n}^{\mathrm{BR}}(nm)
 \ge\E\bigl[\1_{E_t}\{C_t^{\mathrm{pot}}
   +n\mathcal W_{t+1}(M^+/n;\mathbf c_{\mathrm{new}})
   -A_{t+1}\sqrt n-D_{t+1}n^{1-q}\}\mid\mathcal F_t\bigr].
\end{align*}
The event is determined by the completed epoch in the potential
construction. On $E_t$, it is the same observed history as in the
actual system, so conditioning on the next history is legitimate.

To remove the indicator, observe first that
$C_t^{\mathrm{pot}}\le N_A(nm)+\sum_{s\in\cSI,e\in\cE}Y_{s,e}\le Cn$
by the quotas. Also, $m^+$ has bounded total mass, so
$\mathcal W_{t+1}(m^+;\mathbf c_{\mathrm{old}})\le C$: each
session completes at most one turn per epoch, and the total fluid
source mass is bounded. Thus \eqref{eq:br-tracking-errors} yields
\begin{align*}
 \E[\1_{E_t}C_t^{\mathrm{pot}}\mid\mathcal F_t]
 &\ge ng_t-C\sqrt n-Cn\Pp(E_t^c\mid\mathcal F_t),\\
 \E[\1_{E_t}n\mathcal W_{t+1}(M^+/n;\mathbf c_{\mathrm{new}})
       \mid\mathcal F_t]
 &\ge n\mathcal W_{t+1}(m^+;\mathbf c_{\mathrm{old}})
       -CL\sqrt n-Cn\Pp(E_t^c\mid\mathcal F_t).
\end{align*}
Replacing the probability-weighted induction errors by their full
amounts weakens the lower bound. Using
$\Pp(E_t^c\mid\mathcal F_t)\le Cn^{-q}$ and
\eqref{eq:br-fluid-telescope}, we obtain
\[
 V_{t,n}^{\mathrm{BR}}(nm)
 \ge n\mathcal W_t(m;\mathbf c_t^{-,n;t})
      -(A_{t+1}+C')\sqrt n-(D_{t+1}+C'')n^{1-q}.
\]
Set $A_t=A_{t+1}+C'$ and $D_t=D_{t+1}+C''$ to complete the induction.
\end{proof}

\begin{proof}[Proof of Theorem~\ref{thm:br-regret}]
The routing decisions use current arrivals and previously filled
quotas; retention uses the observed boundary states. The emergency
mechanism is causal and enforces capacities. Thus the policy is
nonanticipative and feasible with that mechanism. Taking $q>2$ in
Lemma~\ref{lem:br-rare} gives
$\Pp(\mathsf E_n)\le Cn^{-q}\le n^{-2}$ for all sufficiently large $n$.

The nominal and buffered initial LPs are both feasible, and
\[
 \|\mathbf c_t-\mathbf c_t^{-,n;t}\|_1
 =2E\sum_{\tau=t}^{H-1}(\tau-t+1)\delta_n
 =E(H-t)(H-t+1)b\sqrt{\frac{\log n}{n}}.
\]
Lemma~\ref{lem:br-lipschitz} therefore implies
\[
 n[\mathcal W_t(m;\mathbf c_t)
      -\mathcal W_t(m;\mathbf c_t^{-,n;t})]
 \le C\sqrt{n\log n}.
\]
Combining this reserve loss with Lemma~\ref{lem:br-tracking} gives
\[
 n\mathcal W_t(m)-V_{t,n}^{\mathrm{BR}}(nm)
 \le C\sqrt{n\log n}+C\sqrt n+Cn^{1-q}
 \le C'\sqrt{n\log n}.
\]
These three terms respectively account for the reserve, quota and
transition errors, and exceptional epochs. Since $q>2$, the last
term is $O(1/n)$. Lemma~\ref{lem:br-upper} supplies the comparison
with feasible policies without service interruptions. \end{proof}

\paragraph{Proof of Remark~\ref{rem:br-optimality}.}
Let $C^\pi$ and $D^\pi$ count completed and interrupted turns over
the remaining horizon under policy $\pi$. An interruption ends the
current turn and releases its active memory, and incurs a penalty $\varpi$ independent of $n$, whose value will be specified later. Work already performed
counts against capacity. The affected session may depart or become
idle.

Form the normalized expected flow vector $x\ge 0$ as in
Lemma~\ref{lem:br-upper}, using actual states and decisions and
taking the resource slacks to be expected unused capacities divided
by $n$. Write $Ax=h$ for the nominal fluid equations in the
standard form of Lemma~\ref{lem:br-lipschitz}. Then, for a constant
$k$ depending only on the per-turn resource bounds,
\[
  \|Ax-h\|_1 \le k\,\E D^\pi/n.
\]
To see this, assign each initially active or newly admitted turn
its potential uninterrupted outcome for the current epoch, using
its conditional law at epoch start or at routing, respectively.
These potential outcomes have the means used in the fluid
equations. This calculation allows the subsequent interruption
decision to depend on observed usage. Each interruption changes
one session's successor state, a discrepancy of at most two in
$\ell_1$, and saves at most one turn's bounded execution and
active-memory requirements for that epoch. Arrival and retention
balances hold exactly for the observed states. Summing over epochs
gives the displayed bound.

Let $f$ be the fluid objective vector and $v$ an optimal nominal
dual vertex. By Lemma~\ref{lem:br-lipschitz}, we may take
$\|v\|_\infty\le U$ uniformly over starting epochs and feasible
instances. Since $A^\top v\ge f$ and $x\ge0$,
\[
  f^\top x
  \le v^\top h+v^\top(Ax-h)
  \le \mathcal V_t(m;\mathbf c_t)+Uk\,\E D^\pi/n.
\]
Each interrupted turn must also be subtracted when converting
admissions into completions. The completion accounting identity
therefore gives
\[
  \E[C^\pi-\varpi D^\pi]
  = n[N_A(m)+f^\top x]-(1+\varpi)\E D^\pi
  \le n\mathcal W_t(m)+(Uk-1-\varpi)\E D^\pi.
\]
Thus any fixed $\varpi\ge\max\{0,Uk-1\}$ suffices for the fluid
upper bound.

Finally, choose $b$ large enough to take $q=4$ in
Lemma~\ref{lem:br-rare}. Over the fixed horizon, each session
submits at most one prompt per epoch, so $D^{\rm BR}$ is at most
$(H+1)$ times the total number of initial and newly arriving
sessions. The bounded normalized initial population and
Assumption~\ref{ass:br-scaling} give
$\E[(D^{\rm BR})^2]=O(n^2)$. Since $D^{\rm BR}=0$ outside the
emergency event $\mathsf E_n$, Cauchy--Schwarz gives
\[
  \E D^{\rm BR}
  \le \bigl(\E[(D^{\rm BR})^2]\Pp(\mathsf E_n)\bigr)^{1/2}
  = O(1/n).
\]
The expected penalty is therefore $O(1/n)$ for fixed $\varpi$.
Combining this estimate with \eqref{eq:br-regret} and the fluid
upper bound proves the claim.

\section{Experimental Details and Hyperparameters}
\label{app:exp-details}

\subsection{Telemetry, Prefix-Reuse Measurement, and Baseline Implementation}

\paragraph{Control loop and telemetry.}
In all cluster runs, vLLM asynchronously publishes per-replica running and waiting request counts, pending prefill tokens, KV block utilization, and prefix-reuse counters to the request-scheduler over gRPC at a nominal one-second cadence. Under \textsc{PARK}, replicas additionally send sequenced admission and completion acknowledgments to reconcile in-flight reservations. An independent monitoring process (Prometheus) polls each replica's cumulative generated-token counters every five seconds using actual request timestamps to compute fleet output throughput.

\paragraph{Intrinsic prefix-hit estimator.}
Because vLLM's native prefix hit rate divides reused tokens by total queried prompt tokens, new prompt text registers as a miss even when the entire prior history remains cached. To provide feedback on retained history without external KV transfer, each replica maintains an \emph{intrinsic-hit} estimator using rotating Bloom filters over chained block hashes (excluding preempted lookups). For an admitted lookup $j$, let $\widehat I_j$ denote the reusable prefix length (in blocks) present in the engine's bounded local history and $A_j$ the actual local hit length, yielding recomputation gap $G_j^{\mathrm{rec}} = (\widehat I_j - A_j)_+$. When positive intrinsic mass is observed, the background publisher increments cumulative feedback counters by
\[
 \delta Q_{\mathrm{RS}} = b\sum_j \widehat I_j,
 \qquad
 \delta H_{\mathrm{RS}} = b\sum_j (\widehat I_j - G_j^{\mathrm{rec}}),
\]
where $b$ is the block-to-token scale. The window ratio $\Delta H_{\mathrm{RS}} / \Delta Q_{\mathrm{RS}}$ measures the fraction of locally seen historical prefixes successfully retained, falling back to ordinary hit/query counts when intrinsic mass is zero.

\paragraph{Baseline tie-breaking, staleness, and cohort replenishment.}
All baseline routing controllers and the \emph{SGLang routing} port exclude replicas whose latest heartbeat is older than 5~seconds, holding requests for a 500-ms recheck if no fresh replica is available. Equal routing scores are broken by a rotating cursor over a fixed replica ordering (or lowest replica index in the SGLang port). Upon dispatch, each baseline immediately increments $W_e$ by one and $B_e$ by $U_e$ (or $P$ for \emph{Least prefill}) until the next heartbeat overwrites these estimates; TB's \emph{Shortest queue} wrapper instead retains all dispatches from the preceding one-second window across heartbeats. For the \emph{SGLang routing} port (revision \texttt{a176ba2f7b74}), load is $L_e = R_e + W_e$, complete prompt blocks are $Q = \lfloor P/b \rfloor$, and matched blocks are $h_e = \min\{Q, \lfloor H_e/b \rfloor\}$; the controller falls back to minimum load when $L_{\max} - L_{\min} > 32$ and $L_{\max} > 1.1 L_{\min}$, or when $\max_e h_e / Q \le 0.5$ (evaluated in 32-bit floating point).
All runs in Table~\ref{tab:exp-routing} continuously replenish completed sessions throughout the three-hour window, beginning at the first vLLM request admission with no excluded warmup or concurrency ramp.

{
\subsection{\textsc{PARK} Implementation Details and Hyperparameters}
\label{app:park-implementation}

\paragraph{Forecasts available before dispatch.}
Let $P$ and $O$ denote the realized prompt and output token lengths of a completed turn. Between two consecutive turns, newly appended user or tool-output tokens are measured by $G=(P_{\rm next}-P_{\rm prev}-O_{\rm prev})_+$. For $X\in\{O,G\}$, let $\mathcal D_X$ denote the set of completed model-wide observations and $\mathcal D_{iX}$ the subset belonging to request $i$'s session. Given an initial prior mean $X_0$ and prior weight $\kappa$, the empirical Bayes forecasts are computed via hierarchical shrinkage:
\begin{equation*}
 \overline X=
 \frac{\kappa X_0+\sum_{x\in\mathcal D_X}x}{\kappa+|\mathcal D_X|},
 \qquad
 \widehat X_i=
 \frac{\kappa\overline X+\sum_{x\in\mathcal D_{iX}}x}
      {\kappa+|\mathcal D_{iX}|}.
 \label{eq:online-dual-shrinkage}
\end{equation*}
For a returning session, the next prompt length is forecast as $\widehat P_i=P_i^{\rm prev}+O_i^{\rm prev}+\widehat G_i$; for a new session, $P_i^{\rm prev}=0$ and $\widehat P_i$ equals the smoothed model-wide mean of first-turn prompt lengths, initialized at $P_0$. Both prompt and output forecasts are floored at one token. Let $\widehat h$ denote the smoothed model-wide prefix-reuse fraction, clipped to $[0,1]$, where each completed returning turn contributes its reused prompt tokens divided by $\min\{P_{\rm next},P_{\rm prev}\}$. Setting $o_i=0$ when session $i$ has no current cache owner, the predicted number of reusable tokens at endpoint $e$ is
\begin{equation*}
 \widehat H_{ie}=\mathbf{1}\{e=o_i\}\widehat h
                  \min\{\widehat P_i,P_i^{\rm prev}\}.
 \label{eq:online-dual-reuse}
\end{equation*}
The controller relies on this historical reuse rate because it does not maintain a real-time inventory of resident prefix blocks at each endpoint.

\paragraph{Resource charges and residual capacity.}
Let $C_e^W=\nu_W T_W$ denote the work-backlog budget given nominal processing rate $\nu_W$ and target backlog horizon $T_W$, and let $C_e^M=\rho K_e$ denote the active-memory budget given usable KV-block capacity $K_e$ and safety fraction $\rho$. For a candidate route to endpoint $e$, the estimated resource charges are
\begin{equation*}
 w_{ie}=\min\{(\widehat P_i-\widehat H_{ie})_++\widehat O_i,C_e^W\},
 \qquad m_{ie}=\operatorname{blk}_e(\widehat P_i+\widehat O_i),
 \label{eq:online-dual-costs}
\end{equation*}
where $\operatorname{blk}_e(\cdot)$ converts token counts into physical KV blocks. Capping $w_{ie}$ at $C_e^W$ ensures that an unusually long prompt can still be admitted by utilizing an idle endpoint's entire work budget rather than being permanently blocked. For each resource $j\in\{W,M\}$, let $C_{ej}=C_e^j$ denote capacity, $I_{ej}$ inherited load from active work, and $A_{ej}$ capacity committed to live, unconsumed offers. Residual available capacity is then
\begin{equation*}
 r_{ej}=[C_{ej}-I_{ej}-A_{ej}]_+ .
 \label{eq:online-dual-residual}
\end{equation*}
Inherited work $I_{eW}$ aggregates pending prefill tokens, dispatch reservations awaiting engine acknowledgment, and remaining output estimates for active requests, with sequenced events reconciling dispatches against periodic engine reports. Inherited active memory $I_{eM}$ equals the maximum of active memory reservations and reported occupied blocks after excluding reclaimable idle caches.

\paragraph{Bid and migration safeguards.}
In~\eqref{eq:exp-online-dual-score}, the waiting bonus is $\omega_i=\min\{a_i/\tau_{\rm age},b_{\max}\}$, where $a_i$ denotes the elapsed waiting time (in seconds) since ready request $i$ entered the admission queue. For a migration to an endpoint $e\ne o_i$, the continuation adjustment is given by
\begin{equation*}
 \begin{aligned}
 \Delta_{ie}&=\alpha(F_{io_i}-F_{ie}),\\
 F_{ie}&=\lambda_e\min\{\widehat G_i+\widehat O_i,C_e^W\}
       +\gamma_e\operatorname{blk}_e
           (\widehat P_i+\widehat G_i+2\widehat O_i).
 \end{aligned}
 \label{eq:online-dual-continuation}
\end{equation*}
Only endpoints with status snapshots at most $\tau_{\rm stale}$ seconds old are eligible for routing. While a new session may be routed to any fresh endpoint, a returning session may always use a fresh cache owner $o_i$ and is permitted to migrate to $e\ne o_i$ only if
\begin{equation*}
 \begin{aligned}
 \ell_{o_i}-\ell_e&\ge g_{\rm load}\quad
                       \text{continuously for at least }\tau_{\rm gap},\\
 s_{ie}-s_{io_i}&\ge g_{\rm score},\qquad
 \theta-\theta_i^{\rm bind}\ge\tau_{\rm bind},
 \end{aligned}
 \label{eq:online-dual-migration}
\end{equation*}
where $\ell_e=\max_j I_{ej}/C_{ej}$ is the normalized load at endpoint $e$, $\theta$ is the current time, and $\theta_i^{\rm bind}$ is the session's most recent endpoint-binding time. The load-gap timer resets whenever the load advantage falls below $g_{\rm load}$, a status snapshot becomes stale, or an engine restarts.

\paragraph{One control round.}
Every $\delta$ seconds, the controller reconciles status snapshots with sequenced events, releases completed reservations, and purges expired offers. It excludes stale endpoints and resets an endpoint's dual prices and timers upon detecting an engine restart; if no endpoints are fresh, all pending requests remain on hold. Each session contributes at most one ready request and is eligible only if it has no turn or offer currently in flight. Within each round, the controller executes three steps in sequence:
\begin{enumerate}[leftmargin=*,itemsep=2pt,topsep=3pt]
\item \emph{Measure desired load.}
For each ready candidate $i$ with a nonempty eligible endpoint set $\mathcal E_i$, let $b_i=\max_{e\in\mathcal E_i}s_{ie}$ and let $e_i^\star$ attain this maximum. Candidates are processed in arrival order, breaking ties by lower current maximum normalized desired load at the endpoint and then by fixed endpoint index. Writing $c_{ie,W}=w_{ie}$ and $c_{ie,M}=m_{ie}$, desired load is
\begin{equation*}
 D_{ej}=I_{ej}+A_{ej}
       +\sum_{i:\,\mathcal E_i\ne\varnothing}
          \mathbf{1}\{e_i^\star=e,\, b_i\ge0\}\,c_{ie,j},
 \label{eq:online-dual-demand}
\end{equation*}
where the sum runs over all ready candidates. Note that the nonnegativity filter $b_i\ge0$ applies here even for overdue requests that will receive priority during allocation.
\item \emph{Issue offers.}
Overdue requests ($a_i\ge\tau_{\rm FIFO}$) are served first in arrival order, followed by remaining candidates in decreasing order of $b_i$ (breaking ties by earlier arrival). For each request, eligible endpoints are evaluated in descending order of $s_{ie}$, breaking ties by smaller prospective peak utilization $\max_j(C_{ej}-r_{ej}+c_{ie,j})/C_{ej}$. The controller selects the first unprotected endpoint satisfying $c_{ie,j}\le r_{ej}$ for both resources ($j\in\{W,M\}$) and $s_{ie}\ge0$, waiving only the bid sign test for overdue requests. Once selected, the offer's resource charges are immediately deducted from $r_{ej}$. If an overdue request cannot be accommodated immediately, the controller protects its top-ranked destination against subsequent offers in the same round so that existing work can drain. Unplaced requests remain on hold.
\item \emph{Update prices and reconcile dispatches.}
Initialized at zero, dual prices are updated via~\eqref{eq:exp-online-dual-update} using the desired loads recorded prior to allocation and carried forward to the next round. Issued offers expire after $\tau_{\rm offer}$ seconds. When a client acts on an offer, the controller re-verifies offer validity, snapshot freshness, engine identity, cache ownership, migration safeguards, absence of an in-flight turn, and residual budgets net of other live offers. If all checks pass, the controller binds the request to the endpoint, records a dispatch reservation, and forwards the request; otherwise, it releases the offer and returns the request to the queue.
\end{enumerate}

\paragraph{From a session's retention score to physical blocks.}
The idle age $a$ in~\eqref{eq:bp-retention-index} measures the elapsed time since the session's previous turn completed; when estimating return-time distributions from historical logs, sessions that remain idle at the end of the observation window are treated as right-censored. Substituting the capped return gain and $b_s=\beta L$ into~\eqref{eq:bp-retention-index} yields the length-normalized index
\begin{equation*}
 \widehat I^{\rm proxy}(L,a)
 =\frac{1}{\beta\widehat{\deff}(a,L)}
    \min\left\{\frac{\widehat R}{L},\lambda_e(a_1+a_2L)\right\}.
 \label{eq:bp-retention-proxy}
\end{equation*}
Given a discrete nondecreasing priority map $Q$ with $Q(0)=0$, the aggregate eviction priority of physical KV block $b$ at endpoint $e$ is
\begin{equation*}
 \pi_e(b)=\sum_{i\in\mathcal S_e(b)}Q(\widehat I_i),
 \label{eq:priority_score_block}
\end{equation*}
where $\mathcal S_e(b)$ denotes the set of sessions sharing block $b$ and $\widehat I_i$ is session $i$'s retention index. For instance, if two sessions map to priorities $4$ and $2$, their private blocks receive priorities $4$ and $2$, whereas their shared prefix blocks receive priority $6$. Maintaining unpinned blocks in persistent priority buckets allows the local allocator to reclaim the lowest-priority eligible block in $O(1)$ bucket operations without scanning all idle sessions upon every allocation.

A session with a zero retention index contributes zero weight. For a positive index $I$, the map $Q$ rounds $2^{32}I$ down to a power of two and clips the result to $[1,2^{30}]$. Each reclaimable block is assigned to a bucket according to the sum of its sessions’ weights. Bucket 0 holds blocks of zero total weight. The lower bounds for buckets 1 through 31 are $1,2,4,…,2^{30}$, respectively. The allocator reclaims the block that entered the lowest nonempty bucket earliest. A block that changes buckets joins the tail of its new bucket.

The endpoint queues a session for rescoring when a turn finishes or a return arrives. While a session remains idle, its next refresh is due one second after its previous refresh completes. Each scheduler step starts at most four session lookups and updates at most 256 block contributions, subject to a soft time budget of 0.5 ms.

When a session ends, the worker closes the session at the router. At the endpoint, an idle forecast expires 120 seconds after the last turn completes if no return is pending. If a compatible return is already queued, its idle claim ends when the request attaches to the cached prefix. As sessions gain or lose a claim on a block, the endpoint adds or subtracts each affected session’s contribution and updates the bucket using the new total. For example, a block with contributions of 4 and 2 starts in bucket 3. Withdrawing the 4-unit contribution leaves a total of 2 and moves the block to bucket 2. The cached data remain in place until an allocation reclaims the block with zero active references.

\paragraph{Experimental settings.}
Table~\ref{tab:online-dual-parameters} lists the parameter values used for \textsc{PARK} (Algorithm~\ref{alg:online-dual}) across the cluster experiments. In particular, we set \(\widehat R=131{,}072\) as a fixed reference-value cap, calibrated to the linear work of rebuilding a full context at a reference execution price of one.

\begin{table}[ht]
\centering\small
\caption{Experimental hyperparameters for \textsc{PARK} (Algorithm~\ref{alg:online-dual}).
}
\label{tab:online-dual-parameters}
\begin{tabular}{lll}
\toprule
Parameter & Symbol & Value \\
\midrule
Initial priors (prompt, output, growth) & $(P_0, O_0, G_0)$ & $(2{,}048,\, 256,\, 512)$ tokens \\
Empirical Bayes prior weight & $\kappa$ & $2$ \\
Initial cache-hit fraction prior & $\widehat h_0$ & $0.90$ \\
Work-backlog budget & $C_e^W$ & $131{,}072$ token-work units ($4 \times 32{,}768$) \\
Active KV-memory safety fraction & $\rho$ & $0.90$ \\
Control update interval & $\delta$ & $1.0\text{ s}$ \\
Snapshot staleness threshold & $\tau_{\mathrm{stale}}$ & $5.0\text{ s}$ \\
Admission offer lifetime & $\tau_{\mathrm{offer}}$ & $2.0\text{ s}$ \\
Waiting bonus scale and cap & $(\tau_{\mathrm{age}}, b_{\max})$ & $(60\text{ s},\, 2.0)$ \\
Overdue FIFO threshold & $\tau_{\mathrm{FIFO}}$ & $1{,}500\text{ s}$ \\
Migration continuation weight & $\alpha$ & $0.25$ \\
Minimum normalized load gap & $g_{\rm load}$ & $0.25$ \\
Minimum bid improvement & $g_{\rm score}$ & $0.25$ \\
Minimum load-gap duration & $\tau_{\rm gap}$ & $30\text{ s}$ \\
Minimum time since last binding & $\tau_{\rm bind}$ & $300\text{ s}$ \\
Dual subgradient step size & $\eta$ & $1.0$ \\
Retention cost coefficients & $(a_1, a_2)$ & $(1.0,\, 10^{-6})$ \\
Retention value cap & $\widehat R$ & $131{,}072$ \\
\bottomrule
\end{tabular}
\end{table}
}

\end{document}